\documentclass[twocolumn]{IEEEtran}

\IEEEoverridecommandlockouts    
\usepackage{psfrag}
\usepackage{amsmath,amssymb,amsfonts,bbm}
\usepackage{latexsym}
\usepackage{graphicx}

\usepackage{balance}

\usepackage{colortbl}
\usepackage{fancyhdr}
\usepackage{epstopdf}
\usepackage{float}
\usepackage{subcaption}
\usepackage{hyperref}
\usepackage{booktabs}

\usepackage[usenames,dvipsnames,svgnames]{xcolor}

\newcommand{\proj}{\mathrm{proj}}
\newcommand{\dist}{\mathrm{dist}}
\newcommand{\rge}{\mathrm{rge}}
\newcommand{\Id}{\mathrm{Id}}
\newcommand{\vect}{\mathrm{vec}}
\newcommand{\diag}{\mathrm{diag}}

\newtheorem{theorem}{Theorem}
\newtheorem{definition}{Definition}
\newtheorem{proposition}{Proposition}
\newtheorem{lemma}{Lemma}
\newtheorem{corollary}{Corollary}

\newtheorem{remark}{Remark}

\newtheorem{standing}{Standing Assumption}
\newtheorem{design}{Design choice}

\IEEEoverridecommandlockouts

\usepackage{soul}
\usepackage[ruled]{algorithm2e}

\newcommand{\R}{\mathbb{R}}

\newcommand{\N}{\mathbb{N}}

\newcommand{\mc}{\mathcal}

\newcommand{\bbS}{\mathbb{S}}

\newcommand{\argmin}{\arg \min}
\newcommand{\bs}{\boldsymbol}

\newcommand{\norm}[1]{\left\|#1\right\|}

\begin{document}

\title{ 
Dynamic Control of Agents playing Aggregative Games with Coupling Constraints
}

\author{Sergio Grammatico
\thanks{S. Grammatico is with the Control Systems group, Department of Electrical Engineering, Eindhoven University of Technology, The Netherlands. 
E-mail address: \texttt{s.grammatico@tue.nl}.
}
}
\maketitle

\begin{abstract}
We address the problem to control a population of noncooperative heterogeneous agents, each with convex cost function depending on the average population state, and all sharing a convex constraint, towards an aggregative equilibrium. We assume an information structure through which a central coordinator has access to the average population state and can broadcast control signals for steering the decentralized optimal responses of the agents. We design a dynamic control law that, based on operator theoretic arguments, ensures global convergence to an equilibrium independently on the problem data, that are the cost functions and the constraints, local and global, of the agents. 
We illustrate the proposed method in two application domains: network congestion control and demand side management.
\end{abstract}

\section{Introduction} \label{sec:intro}

\subsubsection*{Motivation}
The problem to coordinate a population of competitive agents arises in several application domains such as the demand side management in the smart grid \cite{mohsenian-rad:10,Saad2012,deng:yang:chen:asr:chow:14,chen:li:louie:vucetic:14,ye:hu:16}, e.g. for thermostatically controlled loads \cite{Li2016, li:zhang:lian:kalsi:16, grammatico:gentile:parise:lygeros:15} and plug-in electric vehicles \cite{ma:callaway:hiskens:13, parise:colombino:grammatico:lygeros:14, ma:zou:ran:shi:hiskens:16, grammatico:16cdc-pev}, demand response in competitive markets \cite{li:chen:dahleh:15}, congestion control for networks with shared resources \cite{barrera:garcia:15}.

The typical challenge in such coordination problems is that the agents are noncooperative, self-interested, yet coupled together, and have local decision authority that if left uncontrolled can lead to undesired emerging population behavior.
From the control-theoretic perspective, the objective is to design a coordination law for steering the strategies of the agents towards a noncooperative game equilibrium.
\smallskip

\subsubsection*{Related literature}
Whenever the behavior of each agent is affected by some aggregate effect of all the agents, which is a typical feature of the mentioned application domains, rather than by agent-specific one-to-one effects, 
\textit{aggregative games} \cite{kukushkin:04, dubey:haimanko:zapechelnyuk:06, dindos:mezzetti:06, jensen:10} offer the fundamentals to analyze the strategic interactions between each individual agent and the entire population, although in the classic literature the analysis is limited to agents with scalar decision variable.

For large, in fact in the limit of infinite, population size, aggregative game setups have been considered as \textit{deterministic mean field games} among agents with strongly convex \textit{quadratic} cost functions \cite{grammatico:parise:lygeros:15, grammatico:parise:colombino:lygeros:16}.

In this paper, we are interested in \textit{generalized} aggregative games for a population of agents with general \textit{convex} functions, \textit{constrained} vector decision variable, and in addition with convex \textit{coupling} (i.e., \textit{shared}) \textit{constraints}.

Generalized games, that is, games among agents with coupling constraints have been intensively studied in the last decade within the operations research community \cite{facchinei:pang, facchinei:kanzow:07} and the control systems one \cite{pavel:07, yin:shanbhag:mehta:11, kulkarni:shanbhag:12, arslan:demirkol:yuksel:15} in relation with duality theory and variational inequalities.

Assessing the convergence of the dynamic interactions among the noncooperative agents towards an equilibrium is one main challenge that arises in (generalized) games. With this aim, best response dynamics and fictitious play with inertia, i.e., gradient update dynamics, have been analyzed and designed, respectively, both in discrete \cite{uryasev:rubinsten:94, marden:arslan:shamma:09} and continuous time setups \cite{shamma:arslan:05, cortes:martinez:15}. In particular, fictitious play with inertia has been introduced to overcome the non-convergence issue of the best response dynamics \cite{shamma:arslan:05}. The common feature of these methods is that the agents implement sufficiently small gradient-type steps, each along the direction of optimality for their local problem. Thus, the \textit{noncooperative} agents shall agree on the sequence of step sizes and exchange truthful information, e.g. with neighboring players, to update their local descent directions. Several distributed algorithms have been proposed for computing the game equilibria, see 
\cite{parise:gentile:grammatico:lygeros:15,Zhu2016,koshal:nedic:shanbhag:16, salehisadaghiani:pavel:16,paccagnan:16} and the references therein.
\smallskip

\subsubsection*{Originality}
In this paper, we consider aggregative games among noncooperative agents that do not exchange information, nor agree on variables affecting their local behavior, with the other (competing) agents.

\begin{figure}[!]
\begin{center}
\includegraphics[width=1\columnwidth]{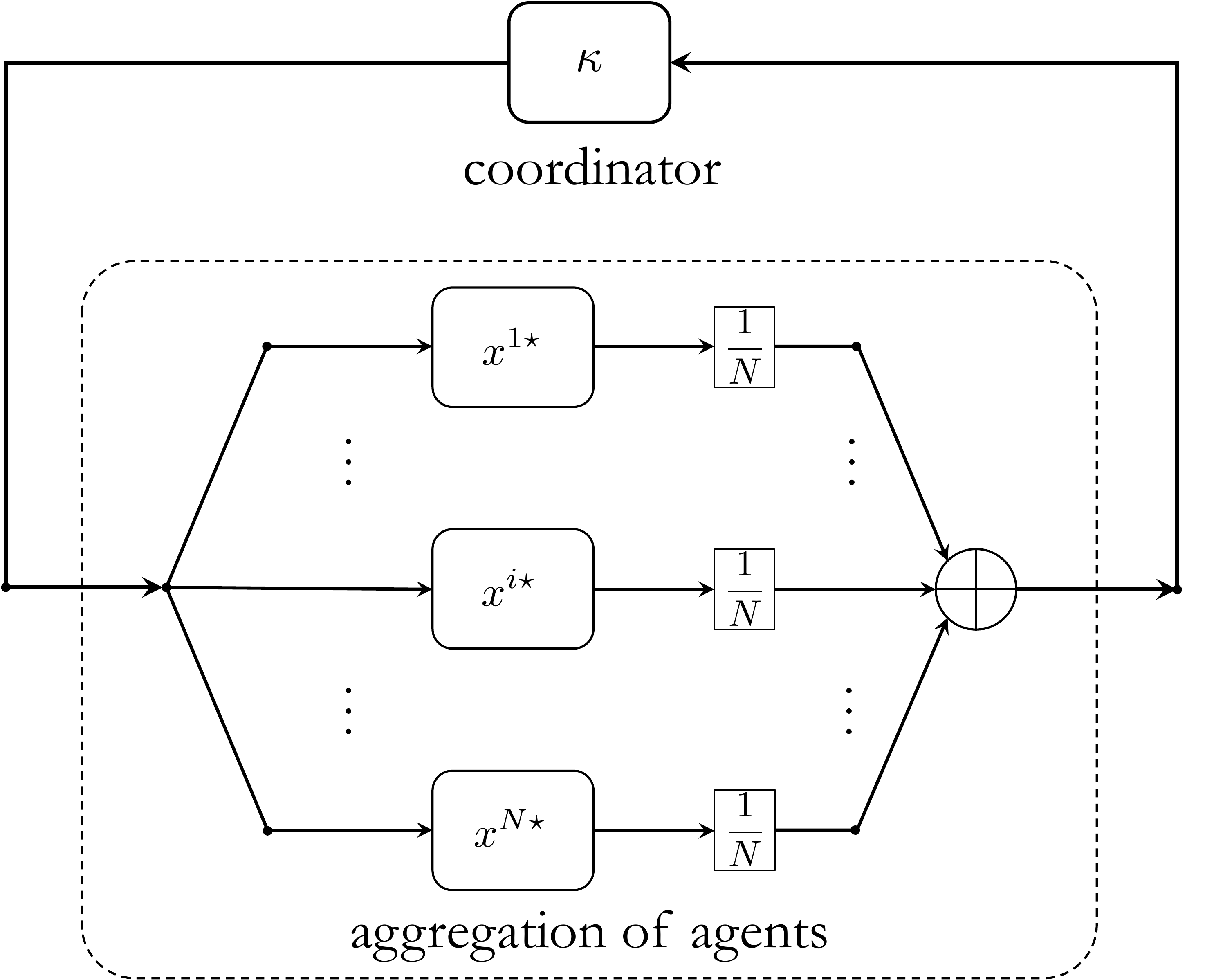}
\end{center}
\caption{Semi-dentralized control architecture. The coordinator $\kappa$ can broadcast to all the agents incentive signals that are designed based on the average among their decentralized optimal responses $(x^{i\star})_{i=1}^{N}$.
\vspace{-0.25cm}
}
\label{fig:scheme}
\end{figure}

Instead, we assume the presence of a central coordinator that controls the decentralized optimal responses of the competitive agents, via the broadcast of incentive signals common to all of them.
Specifically, we design a dynamic control law computing incentives that affect linearly the cost functions of all the agents, simply based on the average among their decentralized optimal responses. The resulting information structure determines the semi-decentralized control architecture illustrated in Figure \ref{fig:scheme}.

Technically, we wish to control the decentralized optimal responses of the agents towards an aggregative equilibrium, that is, a set of agent strategies that are feasible for both the local and the shared constraints, and individually optimal for each agent, given the strategies of all other agents and the control vector associated with the potential violation of the shared constraints.

\smallskip
\subsubsection*{Contribution}
The main contributions and novelties of the paper with respect to the literature are summarized next.
\begin{itemize}
\item We address the general problem to control a population of competitive agents with convex cost functions and constraints coupled together in aggregative form.

\item We discover a nontrivial multivariable mapping with the following two fundamental properties:
\begin{enumerate}
\item[1.] its unique zero is the incentive signal that generates, via the agents' decentralized optimal responses, the desired equilibrium;

\item[2.] there exists a Hilbert space in which the mapping reads as the sum of two monotone operators.
\end{enumerate}

Therefore, splitting methods are applicable for computing the zero of such mapping in a semi-decentralized fashion.

\item We design a dynamic control law with global convergence guarantee for steering the agents' decentralized optimal responses to the desired equilibrium, with minimal information structure, and with no assumption on the problem data, other than convexity.

\item We establish global logarithmic convergence rate under an appropriate selection of the control parameters.

\item We show that our approach is applicable to network congestion control and demand side management.
\end{itemize}

To establish global convergence with minimal information structure, we build upon mathematical tools from variational and convex analysis \cite{rockafellar:wets}, and monotone operator theory \cite{bauschke:combettes}. 

Equilibrium seeking in aggregative games with convex cost functions, convex local constraints and convex coupling constraints has been first studied in \cite{grammatico:16netgcoop}, with \textit{static} control law. In this paper, we enrich the technical setup and study \textit{dynamic} control laws. 
Preliminary versions of some technical results in this paper are in \cite{grammatico:16cdc-convex} where no coupling constraint is considered, and in \cite{grammatico:16cdc-coupling} where the cost functions are assumed to be strongly convex quadratic.

\subsubsection*{Paper organization}
Section \ref{sec:game} define the aggregative game setup. Section \ref{sec:control} presents the novel dynamic control law. The main technical results are shown in Section \ref{sec:main-results} and the designed algorithm is discussed in
Section \ref{sec:discussion}. Section \ref{sec:applicability} illustrates our approach via numerical simulations. Section \ref{sec:conclusion} concludes the paper and points at several research avenues. Some proofs are provided in the Appendix.

\subsection*{Notation}
$\R$, $\R_{> 0}$, $\R_{\geq 0}$ respectively denote the set of real, positive, and non-negative real numbers; $\overline{\R} := \R \cup \{\infty\}$;
$\N$ denotes the set of natural numbers; for $a, b \in \N$, $a \leq b$, $\N[a,b] := [a,b] \cap {\N}$.
$A^\top \in \R^{m \times n}$ denotes the transpose of a matrix $A \in \R^{n \times m}$. Given vectors $x_1, \ldots, x_T \in \R^n$, $\left[ x_1; \ldots; x_T \right] \in \R^{nT}$ denotes $\left[ x_1^{\top}, \cdots, x_T^{\top} \right]^\top \in \R^{n T}$.
Given matrices $A_1, \ldots, A_M$, $\diag\left( A_1, \ldots, A_M\right)$ denotes the block diagonal matrix with $A_1, \ldots, A_M$ in block diagonal positions; given scalars $a_1, \ldots, a_M$, we use the notation $\vect\left( (a_i)_{i=1}^{M} \right) := [a_1, \ldots, a_M]^\top \in \R^M$.
With $\bbS^n$ we denote the set of symmetric $n \times n$ matrices; for a given $Q \in \bbS^n$, {the notations} $Q \succ 0$ ($Q \succcurlyeq 0$) and $Q \in \bbS_{\succ 0}^n$ ($Q \in \bbS_{\succcurlyeq 0}^n$) denote that $Q$ is symmetric and has positive (non-negative) eigenvalues. 
$I$ denotes the identity matrix; $\bs{0}$ ($\bs{1}$) denotes a matrix/vector with all elements equal to $0$ ($1$); to improve clarity, we may add the dimension of these matrices/vectors as subscript. $A \otimes B$ denotes the Kronecker product between matrices $A$ and $B$. Every mentioned set $\mathcal{S} \subseteq \R^n$ is meant to be nonempty. Given $\mathcal{S} \subseteq \R^n$, $A \in \R^{n \times n}$ and $b \in \R^n$, $A \mathcal{S} + b$ denotes the set $\{ A x+b \in \R^n \mid x \in \mathcal{S} \}$; hence 
$\textstyle \frac{1}{N}\sum_{i=1}^{N} \mathcal{S}^i  := \{ \frac{1}{N} \sum_{i=1}^{N} x^i \in \R^n \mid x^i \in \mathcal{S}^i \ \forall i \in \N[1,N] \}$. The notation $\dist\left( x, \mc{S} \right) := \inf_{y \in \mathcal{S}} \left\| x - y\right\|$ denotes the distance of a vector $x \in \R^n$ from a set $\mathcal{S} \subseteq \R^n$.

\subsubsection*{Operator theoretic notations and definitions}
$\mathcal{H}_{Q}$, with $Q \in \bbS_{\succ 0}^n$, denotes the Hilbert space $\R^n$ with inner product $\langle x, y \rangle_{Q} := x^\top Q y$ and induced norm $\left\| x \right\|_Q := \sqrt{ x^\top Q x }$, for all $x,y \in \R^n$; we 
refer to the Hilbert space $\mc{H}_{I}$ whenever not specified otherwise.
Given a function $f:\R^n \rightarrow \overline{\R}$, $\textup{dom}(f) := \{x \in \R^n \mid f(x) < \infty\}$. $f:\R^n \rightarrow \overline{\R}$ is $\ell$-\textit{strongly} \textit{convex}, where $\ell \in \R_{>0}$, if $f(\cdot) - \frac{1}{2} \ell \left\| \cdot \right\|^2$ is convex.
$\text{Id}: \R^n \rightarrow \R^n$ denotes the identity operator.
A mapping $f: \R^n \rightarrow \overline{\R}^n$ is $\ell$-Lipschitz continuous relative to $\mathcal{H}_{Q}$, where $\ell \in \R_{>0}$, if $\left\| f(x) - f(y)\right\|_Q \leq \ell \left\| x-y\right\|_Q$ for all $x,y \in \text{dom}(f)$; $f$ is a contraction (nonexpansive) mapping in $\mathcal{H}_Q$ if it is $\ell$-Lipschitz relative to $\mathcal{H}_{Q}$ with $\ell \in [0,1)$ ($\ell \in [0,1]$).
Given a function $f: \R^n \rightarrow \overline{\R}$, $\partial f: \textup{dom}(f) 
\rightrightarrows {\R}^n$ denotes its subdifferential set-valued mapping \cite{rockafellar:wets}, defined as $\partial f(x) := \{ v \in \R^n \mid f(z) \geq f(x) + v^\top (z-x)  \textup{ for all } z \in \textup{dom}(f) \}$.
A mapping $\mc{T}: \R^n \rightarrow \overline{\R}^n$ is (strictly) \textit{monotone} in $\mathcal{H}_{Q}$ if $\left( \mc{T}(x)-\mc{T}(y)\right)^\top Q \left( x-y\right) \geq 0 \, (> 0)$ for all $x \neq y \in \textup{dom}(\mc{T})$; it is $\ell$-\textit{strongly} \textit{monotone}, where $\ell \in \R_{>0}$, $\mathcal{H}_{Q}$ if $\left( \mc{T}(x)-\mc{T}(y)\right)^\top Q \left( x-y\right) \geq \ell \left\| x-y\right\|_Q^2$ for all $x,y \in \R^n$; 
it is $\beta$-\textit{averaged}, with $\beta \in (0,1)$, if 
$\left\| \mc{T}(x)-\mc{T}(y) \right\|_Q^2 \leq \left\| x-y\right\|_Q^2 
- \frac{1-\beta}{\beta}\left\| \mc{T}(x)-\mc{T}(y)-(x-y) \right\|_Q^2$ for all $x,y \in \R^n$;
it is \textit{firmly nonexpansive} (hence strictly monotone and nonexpansive) if it is $\frac{1}{2}$-averaged; 
it is $\beta$-\textit{cocoercive} (hence strictly monotone), with $\beta \in \R_{>0}$, if the mapping $\beta \, \mc{T}(\cdot)$ is firmly nonexpansive.

\section{Aggregative games with coupling constraints} 
\label{sec:game}

We consider a population of {$N$} agents, where each agent $i \in \N[1,N]$ has strategy (i.e., decision variable) $x^i \in \mc{X}^i \subset \R^n$, and all share the constraint
\begin{equation} \label{eq:shared-constraint}
\textstyle \frac{1}{N} \sum_{i=1}^{N} x^i \in \mc{S},
\end{equation} 
for some set $\mathcal{S} \subseteq \frac{1}{N} \sum_{i=1}^N \mc{X}^i \subset \R^n$.

We assume that each agent $i \in \N[1,N]$ aims at minimizing its local cost function ${J^i}$, which depends on the average among the strategies of all other agents, and in particular  at seeking a strategy $x^i$ such that
\begin{equation} \label{eq:best-response-price}
x^i \in \underset{\, y \in \mc{X}^i}{\argmin} \ \textstyle  J^i\left( y, \, \frac{1}{N} \sum_{j = 1}^N x^j, \, \lambda \right) \,,
\end{equation}
where the argument $\lambda \in \R^{n}$ represents a control vector
that the coordinator agent, introduced later on, can impose on the agents to avoid the violation of the coupling constraint in \eqref{eq:shared-constraint}. 
Equations \eqref{eq:shared-constraint}--\eqref{eq:best-response-price} define a competitive aggregative game. We have an aggregative game since the optimal strategy of each agent depends on the average among the strategies of all agents; the game is competitive aggregative since the cost functions of the agents all depend on a common vector $\lambda$ associated with the coupling constraint in aggregative form. 

Throughout the paper, we assume compactness, convexity and Slater's qualification \cite[\S 5.2.3]{boyd:vandenberghe} of both the individual and the shared constraints, and strong convexity of the cost functions, with linear dependence on the global coupling variable.
Such basic assumptions ensure existence of an equilibrium, and 
that the agents' optimal responses, defined formally in Section \ref{sec:fixed-point}, are single-valued and continuous.

\smallskip
\begin{standing} \label{ass:standing}
\textit{Compactness, convexity, constraint qualification}.
The sets $\{ \mc{X}^i \}_{i=1}^{N}$ and $\mathcal{S} \subseteq \tfrac{1}{N} \sum_{i=1}^{N} \mc{X}^i$ are compact and convex subsets of $\R^n$, and satisfy the Slater's constraint qualification. 
{\hfill $\square$}
\end{standing}
\smallskip

\begin{standing} \label{stand:J}
\textit{Strongly convex cost functions}.
For all $i \in \N[1,N]$, the cost function $J^i: \R^n \times \R^n \rightarrow \overline{\R}$ in \eqref{eq:best-response-price} is defined as
\begin{equation} \label{eq:J}
J^i(y, \sigma, \lambda) := f^i(y) + \left( C \sigma + K \lambda \right)^\top \! y,
\end{equation}
for some function $f^i: \R^n \rightarrow \overline{\R}$ continuous and $\ell$-strongly convex, with $\ell \in \R_{>0}$, $C \in \bbS^{n}$, and invertible $K \in \bbS^{n}$.
{\hfill $\square$}
\end{standing}
\smallskip

In \eqref{eq:J}, the matrix $C$ in \eqref{eq:J} weights the influence of the average among the agents' strategies on each cost function $J^i$, whereas the matrix $K$ in \eqref{eq:J} weights the effect of the vector $\lambda$. In the remainder of the paper, we consider $C$ as part of the given problem data, while $K$ as design choice for the coordinator of the game.

Our goal is to control the strategies of the agents to an aggregative equilibrium, that is, 
a set of strategies and control vector such that: the coupling constraint in \eqref{eq:shared-constraint} is satisfied, and each agent's strategy is optimal given the strategies of all other agents and the control vector.

\smallskip
\begin{definition} \label{def:ANESC}
\textit{Aggregative equilibrium}.
A pair $\textstyle \left( ( \bar{x}^i )_{i=1}^N, \, \bar{\lambda} \right)$ 
is an aggregative equilibrium for the game in \eqref{eq:best-response-price} with coupling constraint in \eqref{eq:shared-constraint} if $\frac{1}{N} \sum_{i=1}^N \bar{x}^i \in \mathcal{S}$, 
for all $i \in \N[1,N]$,
\begin{equation*}
\bar{x}^i \in 
\underset{y \in \mc{X}^i}{\argmin} \ \textstyle J^i\left( y, \, \frac{1}{N} \sum_{j = 1}^{N} \bar{x}^j , \, \bar{\lambda} \right).
\end{equation*}
{\hfill $\square$}
\end{definition}
\smallskip

We formalize next that an aggregative equilibrium exists under the postulated standing assumptions.

\smallskip
\begin{proposition} \label{prop:existence-ANESC}
\textit{Existence of an aggregative equilibrium}.
There exists an aggregative equilibrium for the game in \eqref{eq:best-response-price}  with coupling constraint in \eqref{eq:shared-constraint}.
{\hfill $\square$}
\end{proposition}
\smallskip
\begin{proof}
See Appendix \ref{app:existence-equilibrium}.
\end{proof}

\smallskip
\begin{remark} \label{rem:existence-ANESC}
\textit{Non-uniqueness of aggregative equilibria}.
Uniqueness of the aggregative equilibrium does not necessarily hold. For instance, consider the game with following problem data: $n=1$, $N=2$, $f^1(\cdot) = f^2(\cdot) = \frac{1}{2}\left\| \cdot\right\|^2$, $C = -1$, $\mathcal{X}^1 = \mathcal{X}^2 = \mathcal{S} = [-1, 1]$. The pairs $\left( (1, 1), \, 0\right)$ and $\left( (-1, -1), \, 0\right)$ are aggregative equilibria, independently on the choice of $K$ in \eqref{eq:J}. 
Selecting the best aggregative equilibrium from a global optimization perspective goes beyond the purpose of this paper.
{\hfill $\square$}
\end{remark}
\smallskip

To conclude the section, we note that in the limit of infinite population size, an aggregative equilibrium is a Nash equilibrium with fixed control vector.

\smallskip
\begin{theorem} \label{th:eps-ANESC}
\textit{Aggregative equilibrium versus Nash equilibrium}.
Let the pair $\left( ( \bar{x}^i )_{i=1}^{N}, \, \bar{\lambda} \right)$ be an aggregative equilibrium, and define
\begin{multline*}
 \varepsilon_N := \max_{ i \in \N[1,N] } \,
\dist\left( \, \bar{x}^i \, , \right. \\
\left. \textstyle \underset{y \in \mathcal{X}^i}{\argmin} \, J^i\left( y, \frac{1}{N}\left( y + \sum_{j\neq i}^{N} \bar{x}^j \right), \bar{\lambda} \right) \right) .  
\end{multline*}
Assume that there exists a compact set $\mc{X} \subset \R^n$ such that $\mc{X}^i \subseteq \mc{X}$ for all $i \in \N[1,N]$ and $N \in \N$. Then, there exists $c \in \R_{>0}$ such that $\varepsilon_N \leq c/N$ for all $N \in \N$.
{\hfill $\square$}
\end{theorem}
\smallskip
\begin{proof}
See Appendix \ref{app:theorem-1}.
\end{proof}

\section{Dynamic control of the agents' decentralized optimal responses}
\label{sec:control}

\subsection{Fixed points of the aggregation mapping}
\label{sec:fixed-point}
For seeking an aggregative equilibrium, we assume that an agent $i$ cannot exchange information, {nor has prior knowledge}, on the strategies of all other (competing) agents. 
Instead, we assume that each individual agent responds optimally to incentive signals $u \in \R^n$ according to the information structure in Figure \ref{fig:scheme}.
Formally, for all $i \in \N[1,N]$, we define the agent \textit{optimal response} mapping $x^{i \star}: \R^n \rightarrow \mc{X}^i$ as
\begin{align} \label{eq:optimal-responses}
\begin{split}
x^{i \star }(u) & := \underset{ \, y \in \mc{X}^i }{\argmin } \ f^i(y) +  u^\top y,
\end{split}
\end{align}
and the aggregation mapping $\mc{A}: \R^n \times \R^n \rightarrow \frac{1}{N} \sum_{i=1}^{N} \mc{X}^i$ as 
the average among the optimal responses of agents to the incentive signal $u(\sigma, \lambda) = C \sigma + K \lambda$, i.e., 
\begin{equation} \label{eq:aggregation-mapping}
\textstyle \mc{A}(\sigma, \lambda) := \frac{1}{N} \sum_{i=1}^{N} x^{i \star }(C \sigma + K \lambda).
\end{equation}

Note that if $\bar{\sigma} = \mathcal{A}\left( \bar{\sigma}, \bar{\lambda}\right)$ for some $\bar{\lambda} \in \R^n$, then $\bar{\sigma} = \frac{1}{N} \sum_{i=1}^{N} \bar{x}^i$, 
with shorthand notation $\bar{x}^i := x^{i \star }(C \bar\sigma + K \bar{\lambda})$. It follows immediately from Proposition \ref{prop:existence-ANESC} that such a pair $(\bar{\sigma}, \bar{\lambda} )$ exists; uniqueness depends however on the choice of $K$ as established later  in Proposition \ref{prop:existence-uniqueness}, Section \ref{sec:main-results}.

\begin{figure*}[t]
\begin{align} \label{eq:forward-backward} \tag{$*$}
\begin{split}
& \left[ \begin{matrix} \sigma_{(t+1)} \\ \lambda_{(t+1)} \end{matrix} \right] \ = \ \kappa \left( t, \,  \left[ \begin{matrix} \sigma_{(t)} \\ \lambda_{(t)} \end{matrix} \right]\right) \ := \
\,  
(1-\alpha_t) \left[ \begin{matrix} \sigma_{(t)} \\ \lambda_{(t)} \end{matrix} \right] + \alpha_t \, \left( I + \epsilon M  \right)^{-1} \left( \left[ \begin{matrix} \sigma_{(t)} \\ \lambda_{(t)} \end{matrix} \right] - \epsilon \, \Gamma \left( \left[ \begin{matrix} \sigma_{(t)} \\ \lambda_{(t)} \end{matrix} \right]\right)   \right)
\end{split}
\end{align}
\hrule
\end{figure*}

Therefore, if $\bar{\sigma} = \mathcal{A}\left( \bar{\sigma}, \bar{\lambda}\right)$, then the pair $\left( ( \bar{x}^i )_{i=1}^N, \, \bar{\lambda}\right)$ is in fact an aggregative equilibrium. It follows that we can control the agents' optimal responses, e.g. via dynamic updates of {their argument}, to a set of strategies whose average is a fixed point of the aggregation mapping (with respect to the first argument) within the coupling constraint set.

\subsection{From fixed points to zeros}
Informally speaking, the objective is to find a pair $(\bar \sigma, \bar\lambda)$ such that $\, \bar\sigma = \mc{A}(\bar{\sigma}, \bar{\lambda}) = x^0$, for some $x^0\in \mc{S}$. 
Since $\mathcal{A}$ depends on two arguments, it follows naturally that $x^0$ is designed as a mapping that depends on the same arguments.
With this aim, we translate the problem into that of finding a zero of an appropriate multivariable mapping via semi-decentralized iterations.

Among all possible design choices, let us {define} the mapping $x^{0 \star}: \R^n \times \R^n \rightarrow \mc{S}$ as
\begin{equation} \label{eq:optimal-response-0}
x^{0 \star }(\sigma, \lambda) := \underset{ y \in \mathcal{S} }{\argmin} \ 
\textstyle \frac{1}{2} y^\top y + \, \left( K(\sigma-\lambda) \right)^\top \! y.
\end{equation}

Remarkably, we notice that a pair $\left( \bar{\sigma}, \bar{\lambda} \right)$ 
satisfies $\bar\sigma = \mathcal{A}\left(\bar\sigma, \bar\lambda\right) = x^{0 \star}(\bar\sigma, \bar\lambda)\in \mathcal{S}$ if $\left[ \bar{\sigma} \, ; \, \bar{\lambda} \right]$ is a zero of the mapping $\Theta: \R^{2n} \rightarrow \R^{2n}$ defined as
\begin{align} \label{eq:Theta}
\begin{split}
\Theta\left( \left[ \begin{matrix} \sigma \\ \lambda \end{matrix} \right] \right)  & :=
 \left[ \begin{matrix}  \sigma - \mc{A}(\sigma, K \lambda) \\ \sigma -2\mc{A}(\sigma, K \lambda) + x^{0 \star }(\sigma, \lambda) \end{matrix} \right] \\
\medskip 
  & \, =  \left[ 
\begin{matrix}
I_n & 0 \\
I_n & 0
\end{matrix}
\right] \left[ \begin{matrix}
\sigma \\
\lambda
\end{matrix}
\right] -
\left[ \begin{matrix}  \mc{A}(\sigma, K \lambda) \\  2\mc{A}(\sigma, K \lambda) - x^{0 \star }(\sigma, \lambda) \end{matrix} \right] \smallskip \\
  & =:
\left( M + \Gamma \right) \left( \left[ \begin{matrix} \sigma \\ \lambda \end{matrix} \right] \right) \,,
\end{split}
\end{align}
where the matrix gain $K$ is to be chosen, and we defined the matrix $M \in \R^{2n \times 2n}$ and the mapping $\Gamma: \R^{2n} \rightarrow \R^{2n}$ as
\begin{align} \label{eq:M-in-Theta}
M & := \left[ \begin{matrix}
I_n & 0 \\
I_n & 0
\end{matrix}
\right], \\
\Gamma\left( \left[ \begin{matrix} \sigma \\ \lambda \end{matrix} \right] \right) & :=
 -\left[ \begin{matrix}  \mc{A}(\sigma, K \lambda) \\  2\mc{A}(\sigma, K \lambda) - x^{0 \star }(\sigma, \lambda) \end{matrix} \right].
 \label{eq:Gamma}
\end{align}

\subsection{Dynamic control as zero finding algorithm}

In general, computing a zero of a multivariable nonlinear mapping such as $\Theta = M + \Gamma$ in \eqref{eq:Theta} is a challenging task. However, for the sum of monotone mappings there exist iterative algorithms with global convergence guarantee \cite[Chapter 25]{bauschke:combettes}. Inspired by the forward-backward algorithm \cite[Equation 25.26]{bauschke:combettes}, we propose the dynamic control law $\kappa$ in \eqref{eq:forward-backward} for computing a zero of $\Theta = M + \Gamma$ in \eqref{eq:Theta}, where $\epsilon >0$ is sufficiently small and the averaging step sizes $\left( \alpha_t \right)_{t=0}^{\infty}$ are chosen as follows.
\smallskip
\begin{design} \label{des:alpha}
The sequence $\left( \alpha_t \right)_{t=0}^{\infty}$ in \eqref{eq:forward-backward} is such that $\alpha_t \in (0, 3/2)$ for all $t \in \N$ and $\sum_{t=0}^{\infty} \alpha_t \left( \frac{3}{2} - \alpha_t \right) = \infty$.
\hfill $\square$
\end{design}
\smallskip

{Suitable choices for the sequence $(\alpha_t)_{t \in \N}$ that satisfy the design condition stated above are
$\alpha_t = 1$ and $\alpha_t = 1/(t+1)$, for all $t \in \N$.}
The {proposed} dynamic control scheme is summarized in Algorithm \ref{alg:aggregative-control}.

\smallskip
\begin{algorithm}\caption{Dynamic control of the decentralized optimal responses by the competitive agents.} \label{alg:aggregative-control}

\smallskip
{Initialization}: $t \leftarrow 0$; 

\smallskip
$\bullet \,$ The coordinator chooses $( \sigma_{(0)}, \, \lambda_{(0)} ) \in \mc{S} \times \R^n$.

\medskip
{Iterate until convergence}: 

\smallskip
$\bullet \,$ The coordinator broadcasts
{$$ u_{(t)} := C \sigma_{(t)} +  K \lambda_{(t)} $$}
to all agents, and computes $x^{0 \star}\left( \sigma_{(t)}, \lambda_{(t)} \right)$ from \eqref{eq:optimal-response-0}.

\smallskip
$\circ \,$ The agents compute in parallel 
{$x^{i \star}\left( u_{(t)} \right)$} 
from \eqref{eq:optimal-responses}, for all $i \in \N[1,N]$.

\smallskip
$\bullet \,$ The coordinator receives 
{$\mc{A}\left( \sigma_{(t)}, \lambda_{(t)} \right)$} 
as in \eqref{eq:aggregation-mapping}, computes 
$\Gamma\left( [\sigma_{(t)} \, ; \, \lambda_{(t)} ]\right)$ from \eqref{eq:Gamma}, and from \eqref{eq:forward-backward}
$$\left[ \sigma_{(t+1)} \, ; \,  \lambda_{(t+1)} \right] = 
\kappa \left( t, \,  \left[ \sigma_{(t)} \, ; \,  \lambda_{(t)} \right]\right).$$

\smallskip
$t \leftarrow t+1$.
\end{algorithm}

\section{Global convergence} \label{sec:main-results}

\subsection{Statement of the main results}

The mapping $\Theta$ in \eqref{eq:Theta} reads as the sum of the linear, hence continuous, mapping {$M \cdot$} and the mapping {$\Gamma(\cdot)$} in \eqref{eq:Gamma}. With the aim of applying \cite[Theorem 25.8]{bauschke:combettes}, in the following we show that by choosing the matrix gain $K$ in \eqref{eq:optimal-responses} appropriately, $M$ is monotone and $\Gamma$ is $\beta$-cocoercive, that is, $\beta \, \Gamma(\cdot)$ is firmly nonexpansive, in some Hilbert space. Consequently, we derive a dynamic control law that ensures global convergence of the controlled decentralized optimal responses to a set of strategies whose average is a fixed point of the aggregation mapping with respect to its first argument.

\smallskip
\begin{design} \label{des:K}
{The matrix $K$} in \eqref{eq:Theta} is chosen such that $K \succ 0$ and $C+K \succ 0$.
\hfill $\square$
\end{design}
\smallskip

\smallskip
\begin{theorem} \label{th:monotone}
\textit{Monotonicity}.
Under design choice \ref{des:K}, the linear mapping $M$ in \eqref{eq:M-in-Theta} is monotone in $\mathcal{H}_{P}$, and the mapping $\Gamma$ in \eqref{eq:Gamma} is $\beta$-cocoercive, 
{hence strictly monotone,} in $\mathcal{H}_{P}$, where
\begin{equation} \label{eq:P-C-K}
P := \left[ \begin{matrix} C+2K & -K \\ -K & K \end{matrix}\right] \succ 0, \quad \beta := \frac{\ell }{6 \norm{P}} > 0.
\end{equation}
{$\hfill \square$}
\end{theorem}
\smallskip
\begin{proof}
See Section \ref{sec:proof-convergence}.
\end{proof}
\smallskip

\smallskip
\begin{proposition} \label{prop:existence-uniqueness}
\textit{Existence and uniqueness}. Under design choice \ref{des:K}, 
$\exists ! \,  \left[ \bar\sigma \, ; \bar\lambda \right] \in \textup{zer}\left(\Theta\right)$,
with $\Theta$ as in \eqref{eq:Theta}.
{\hfill $\square$}
\end{proposition}

\smallskip
\begin{proof}
Existence follows immediately from Proposition \ref{prop:existence-ANESC}. The mapping $\Theta$ is the sum of monotone and strictly monotone mappings by Theorem \ref{th:monotone}, hence {it is} strictly monotone in $\mc{H}_P$ \cite[Exercise 12.4 (c)]{rockafellar:wets}, with $P$ in \eqref{eq:P-C-K}, hence uniqueness holds \cite[Proposition 23.35]{bauschke:combettes}.
\end{proof}
\smallskip

\smallskip
\begin{design} \label{des:epsilon}
The constant $\epsilon$ in \eqref{eq:forward-backward} is such that $\epsilon \in (0, \beta)$, with $\beta$ in \eqref{eq:P-C-K}.
\hfill $\square$
\end{design}
\smallskip

\smallskip
\begin{theorem} \label{th:forward-backward}
\textit{Global convergence}.
Under design choices \ref{des:alpha}--\ref{des:epsilon}, 
the sequence 
$\left( \left[ \sigma_{(t)} \, ; \, \lambda_{(t)}\right]\right)_{t=0}^{\infty}$
defined in \eqref{eq:forward-backward}
converges, for any initial condition, to {the} zero of $\Theta$ in \eqref{eq:Theta}, with $\mathcal{A}$ as in \eqref{eq:aggregation-mapping} and $x^{i \star}$ as in \eqref{eq:optimal-responses} for all $i \in \N[1,N]$.
{$\hfill \square$}
\end{theorem}
\smallskip

\begin{proof}
The assumptions of \cite[Theorem 25.8 ($A:=M \Rightarrow J_{ \epsilon A} = (I + \epsilon M )^{-1}$)]{bauschke:combettes} are verified as follows. $M$ is continuous and monotone in $\mc{H}_P$ due to Theorem \ref{th:monotone}, hence maximally monotone \cite[Corollary 20.25]{bauschke:combettes}. $\Gamma(\cdot)$ in \eqref{eq:Gamma} is $\beta$-cocoercive in $\mc{H}_P$ according to Theorem \ref{th:monotone}. The sequence $(\alpha_t)_{t=0}^{\infty}$ is chosen as in \cite[Theorem 25.8 ($\gamma \leq \beta $)]{bauschke:combettes}, and the existence of a zero of $\Theta$ holds by Proposition \ref{prop:existence-uniqueness}.
\end{proof}
\smallskip

We conclude the subsection by quantifying the global convergence rate. Since in general this might depend on the 
chosen sequence $(\alpha_t)_{t \in \N}$, let us focus on the case $\alpha_t = \bar{\alpha} \in (0,1)$ for all $t \in \N$, for which we establish global logarithmic convergence rate.

\smallskip
\begin{design} \label{des:alpha=1}
The sequence $\left( \alpha_t \right)_{t=0}^{\infty}$ in \eqref{eq:forward-backward} is such that $\alpha_t = \bar{\alpha} \in (0,1]$ for all $t \in \N$ .
\hfill $\square$
\end{design}
\smallskip

\smallskip
\begin{theorem} \label{th:convergence-rate}
\textit{Global logarithmic convergence rate}.
Under design choices \ref{des:K}--\ref{des:alpha=1}, 
the sequence $\left( \sigma_{(t)}, \, \lambda_{(t)}\right)_{t=0}^{\infty}$ defined in \eqref{eq:forward-backward} is such that, for all $t \in \N$ and any initial condition, 
\begin{equation} \label{eq:logarithmic-rate}
\left\| 
\left[ 
\begin{smallmatrix}
\sigma_{(t+1)} \\
\lambda_{(t+1)}
\end{smallmatrix}
\right] 
- \left[ 
\begin{smallmatrix}
\sigma_{(t)} \\
\lambda_{(t)}
\end{smallmatrix}
\right] \right\|_{P}^2 
\leq \,  \frac{ \textstyle \frac{3}{\bar{\alpha}} -1 }{t+1} \,
\left\| 
\left[ 
\begin{smallmatrix}
\sigma_{(0)} \\
\lambda_{(0)}
\end{smallmatrix}
\right]
 - \textup{zer}\left( \Theta \right) \right\|_{P}^2,
\end{equation}
where $\Theta$ is as in \eqref{eq:Theta}, $P$ as in \eqref{eq:P-C-K}, $\mathcal{A}$ as in \eqref{eq:aggregation-mapping} and $x^{i \star}$ as in \eqref{eq:optimal-responses} for all $i \in \N[1,N]$.
{$\hfill \square$}
\end{theorem}
\smallskip
\begin{proof}
See Section \ref{sec:proof-linear-rate}.
\end{proof}
\smallskip

\smallskip
\begin{corollary}
\textit{Global convergence to an aggregative equilibrium}.
Under design choices \ref{des:K}--\ref{des:alpha=1}, for any initial condition, the sequence 
$\textstyle \left( (x^{i \star}( C \sigma_{(t)} + K \lambda_{(t)}))_{i=1}^N, \, \lambda_{(t)} \right)_{t=0}^{\infty}$
 defined from \eqref{eq:optimal-responses} and \eqref{eq:forward-backward} converges with the logarithmic rate in \eqref{eq:logarithmic-rate} to an aggregative  equilibrium for the game in \eqref{eq:best-response-price} with coupling constraint in \eqref{eq:shared-constraint}.
{$\hfill \square$}
\end{corollary}
\smallskip
\begin{proof}
It follows immediately from Theorems \ref{th:forward-backward}--\ref{th:convergence-rate}.
\end{proof}

\subsection{Proof of Theorem \ref{th:monotone} (Monotonicity)}
\label{sec:proof-convergence}
First, $M$ in \eqref{eq:M-in-Theta} is monotone in $\mc{H}_P$ as \cite[Lemma 3]{grammatico:parise:colombino:lygeros:16} 
$$ \left[ \begin{matrix} I & \bs{0} \\ I & \bs{0} \end{matrix} \right]^\top \left[ \begin{matrix} C+2K & -K \\ -K & K \end{matrix} \right] + \left[ \begin{matrix} C+2K & -K \\ -K & K \end{matrix} \right] \left[ \begin{matrix} I & \bs{0} \\ I & \bs{0} \end{matrix} \right] \succcurlyeq 0. $$

We proceed with two statements that are exploited later on.

\smallskip
\begin{lemma} \label{lem:strong-convexity-monotonicity}
If a function $f:\R^n \rightarrow \overline{\R}$ is $\ell$-strongly convex, $\ell \in \R_{>0}$, then: 
$\partial f$ is $\ell$-strongly monotone, and $\left( \partial f \right)^{-1}$ is everywhere single-valued, globally $(1/\ell)$-Lipschitz continuous, $\ell$-cocoercive, and strictly monotone.
{$\hfill \square$}
\end{lemma}

\smallskip
\begin{proof}
$\partial f$ is $\ell$-strongly monotone by \cite[Exercise 12.59]{rockafellar:wets}, and equivalently $\left( \partial f \right)^{-1}$ is $\ell$-cocoercive \cite[p. 1021, Equation (18)]{combettes:00}.
$\left( \partial f \right)^{-1}$ is everywhere single-valued, globally $(1/\ell)$-Lipschitz continuous by \cite[Proposition 12.54]{rockafellar:wets}.
Finally, we show that $\left( \partial f \right)^{-1}$ is strictly monotone. For all $\xi, \zeta \in \textup{dom}\left( \partial f \right)$ such that $\xi \neq \zeta$, we have
$(v - w )^\top \left( \xi - \zeta \right) \geq \ell \left\| \xi - \zeta\right\|^2 > 0
$ for all $v \in \partial f( \xi ), \, w \in \partial f( \zeta )$.
In particular, since $\left( \partial f \right)^{-1}$ is everywhere single-valued, for all $x, y \in \textup{rge}\left( \partial f \right) = \textup{dom}( \left( \partial f \right)^{-1} ) $ there exist $\xi = \left( \partial f \right)^{-1}(x)$, $ \zeta = \left( \partial f \right)^{-1}(y)$, such that $x \in \partial f(\xi)$, $y \in \partial f(\zeta)$, and hence
$
\left( x- y \right)^\top \left( \left( \partial f\right)^{-1}( x ) - \left( \partial f\right)^{-1}( y ) \right) \geq \ell \left\| \xi - \zeta\right\|^2 > 0
$.
\end{proof}

\smallskip
\begin{lemma} \label{lem:optimizer}
Let the function $f: \R^n \rightarrow \overline{\R}$ be $\ell$-strongly convex, $\ell \in \R_{>0}$. Then for any $A \in \R^{n \times m}$, the mapping 
\begin{equation} \label{eq:argmin_inverse_subdiff}
x^{\star}(\cdot) := \underset{ y \in \R^n}{ \, \arg\min} \, f(y) + \left( A \, \cdot \right)^\top y  \, = \, (\partial f)^{-1}\left( -A \, \cdot \right)
\end{equation}
is $(\left\| A \right\|/\ell)$-Lipschitz continuous.
{\hfill $\square$}
\end{lemma}
\smallskip

\begin{proof}
By Lemma \ref{lem:strong-convexity-monotonicity} the mapping $\left( \partial f\right)^{-1}$ is $(1/\ell)$-Lipschitz continuous. The affine mapping $- A \, \cdot $ is $\left\| A \right\|$-Lipschitz continuous, hence the composed mapping 
$\left( \partial f\right)^{-1}\left( -A \, \cdot\right)$ is $(\left\| A \right\|/\ell)$-Lipschitz continuous.
Equation \eqref{eq:argmin_inverse_subdiff} follows from the Fermat's rule \cite[Theorem 16.2, Proposition 26.1]{bauschke:combettes}, i.e., $0 \in \partial \left( f( \cdot ) + \left( A \ z \right)^\top \cdot\right)\left(x^{\star}(z)\right) \in \partial f ( x^{\star}(z) ) + A z $, hence $-A z \in \partial f\left( x^{\star}(z) \right)$ for all $z \in \R^n$. The second equation in \eqref{eq:argmin_inverse_subdiff} follows by applying $\left( \partial f \right)^{-1}$ to both sides of the last inclusion.
\end{proof}
\smallskip

It follows from Lemma \ref{lem:optimizer} that, for all $i \in \N[1,N]$, the optimal response from \eqref{eq:optimal-responses} reads as
\begin{equation*}
{x^{i \star}(C \sigma + K \lambda)} = (\partial f^i)^{-1}\left( [-C, \, -K] \left[ \begin{matrix} \sigma \\ \lambda \end{matrix}\right] \right),
\end{equation*}
and analogously, the mapping $x^{0 \star}$ in \eqref{eq:optimal-response-0} reads as
\begin{equation*}
x^{0 \star}( \sigma, \lambda ) = (\partial f^0)^{-1}\left( [-K, \, K] \left[ \begin{matrix} \sigma \\ \lambda \end{matrix}\right] \right),
\end{equation*}
where $f^0(y) := \frac{1}{2} y^\top y + \delta_{\mc{S}}(y)$.

In view of $\Gamma$ in \eqref{eq:Gamma}, for all $i \in \N[1,N]$, let us define the mapping $\Gamma^i: \R^{2n} \rightarrow \R^{2n}$ as
\begin{equation} \label{eq:Gamma-i-mapping}
\begin{array}{l}
\Gamma^i([\sigma \, ; \, \lambda]) := -\left[ \begin{matrix} 
{x^{i \star}(C \sigma + K \lambda)} \\ 
2 {x^{i \star}(C \sigma + K \lambda)} - x^{0 \star}(\sigma,\lambda) \end{matrix} \right] = \smallskip \\
- \left[ \begin{matrix} I_n & \bs{0} \\ 2 I_n & -I_n \end{matrix} \right] 
\left[
\begin{matrix} (\partial f^i)^{-1} \! & \bs{0} \\
 \bs{0} & \! (\partial f^0)^{-1}
\end{matrix}
\right] 
\left( \left[ \begin{matrix} 
-C & -K \\ -K & K \end{matrix} \right] \left[ \begin{matrix} \sigma \\ \lambda \end{matrix} \right] \right) 
\medskip
\end{array}
\end{equation}
so that $\Gamma(\cdot) = \frac{1}{N} \sum_{i=1}^{N} \Gamma^i(\cdot)$.
Note that the mapping $\diag\left( \partial f^i, \, \partial f^0\right)$ is $\gamma$-strongly monotone with $\gamma := \min\{ \ell,1 \}$, thus the mapping $\diag\left( ( \partial f^i )^{-1} , \, ( \partial f^0 )^{-1} \right)$ in \eqref{eq:Gamma-i-mapping} is $\gamma$-cocoercive and $(1/\gamma)$-Lipschitz continuous due to Lemma \ref{lem:strong-convexity-monotonicity} and \cite[Proposition 20.23]{bauschke:combettes}. In the rest of the proof, we exploit the following result, which is a variant of \cite[Proposition 4.5]{bauschke:combettes}.

\smallskip
\begin{lemma} \label{lem:coco-preservation}
Let $\mathcal{M}: \R^m \rightarrow \R^m$ be a $\gamma$-cocoercive mapping, $\gamma \in \R_{>0}$, and $A, B \in\R^{m \times m}$ be invertible matrices. If $A^{-\top} B \in \mathbb{S}_{\succ 0}^{m}$, then the mapping $A \, \mathcal{M}\left( B \, \cdot \right)$ is $\eta$-cocoercive in $\mathcal{H}_{A^{-\top} B}$ with  
$\textstyle \eta := \gamma / ( \norm{A}^2 \norm{A^{-\top} B} )$.
{\hfill $\square$}
\end{lemma}
\smallskip
\begin{proof}
Since $\mc{M}$ is $\gamma$-cocoercive, for all $x, y \in \R^m$:
\begin{align*}
& \left( A \mc{M}( B x ) - A \mc{M}( B y ) \right)^\top A^{-\top} B (x-y) \\
& \, = \left( \mc{M}(Bx) - \mc{M}(By) \right)^\top B (x-y) \\
& \, \geq \gamma \left\| \mc{M}(B x) - \mc{M}(By) \right\|^2 \\
& \, \geq \frac{\gamma}{ \norm{A}^2 } \left\| A \mc{M}(B x) - A \mc{M}(By) \right\|^2 \\
& \, \geq \frac{\gamma}{ \norm{A}^2 \norm{A^{-\top} B} } \left\| A \mc{M}(B x) - A \mc{M}(By) \right\|_{A^{-\top} B }^2.
\end{align*}
\end{proof}
\smallskip

We now apply Lemma \ref{lem:coco-preservation} to the mapping $\Gamma^i(\cdot)$ in \eqref{eq:Gamma-i-mapping}. 
Namely, we consider $m=2n$ and the matrices
\begin{equation*}
A := -\left[ \begin{matrix} I_n & \bs{0}_{n \times n} \\ 2 I_n & -I_n \end{matrix} \right], \ 
{B} := \left[ \begin{matrix} -C & -K \\ -K & K \end{matrix} \right],
\end{equation*}
and derive
\begin{align*}
\begin{split}
P := {A}^{-\top} B &= -\left[ \begin{matrix} 
I & \bs{0} \\ 2 I & -I \end{matrix} \right]^{-\top}  
 \, \left[ \begin{matrix} -C  & -K \\ -K & K \end{matrix} \right] \smallskip\\
& = \left[ \begin{matrix} I & 2I \\ \bs{0} & -I \end{matrix} \right] \,  
\left[ \begin{matrix} C  & K \\ K & -K \end{matrix} \right] \smallskip\\
& = \left[ \begin{matrix} C+2K  & -K \\ -K & K \end{matrix} \right] \smallskip\\
& = \left[ \begin{matrix} C+ (1-\epsilon) K & \bs{0} \\ \bs{0} & \bs{0} \end{matrix} \right]  + \left[ \begin{matrix} 1+\epsilon & -1 \\ -1 & 1\end{matrix} \right] \otimes K,
\end{split}
\end{align*}
where $\epsilon > 0$ is chosen such that $C+ (1-\epsilon)K \succcurlyeq 0$.

Since $\left[ \begin{smallmatrix} 1+\epsilon & -1 \\ -1 & 1\end{smallmatrix} \right] \succ 0$ and $C \in \bbS^n$, $K \succ 0$ and $C + K \succ 0$ ensure that $B$ is invertible and ${A}^{-\top} B \succ 0$. By Lemma \ref{lem:coco-preservation}, this implies that, for all $i \in \N[1,N]$, $\Gamma^i(\cdot)$ is $\beta$-cocoercive in $\mc{H}_{ {A}^{-\top} B}$, where ${A}^{-\top} B = P$ and 
$\ell/( \norm{A}^2 \norm{P} ) = \ell/( (3+2\sqrt{2}) \norm{P} ) \geq \ell/(6 \norm{P}) =: \beta$ in \eqref{eq:P-C-K}. In turn, $\Gamma(\cdot) = \frac{1}{N} \sum_{i=1}^{N} \Gamma^i(\cdot)$ is also $\beta$-cocoercive \cite[Example 4.31]{bauschke:combettes}.

Since all the mappings $\{ \Gamma^i \}_{i=1}^{N}$ are strictly monotone in $\mc{H}_P$, it follows that $\Gamma$ is strictly monotone as well \cite[Exercise 12.4 (c)]{rockafellar:wets}; in fact, for all $i \in \N[1,N]$, 
$\diag\left( (\partial f^i)^{-1}, (\partial f^0)^{-1} \right) $ in \eqref{eq:Gamma-i-mapping} is strictly monotone by Lemma \ref{lem:strong-convexity-monotonicity} and \cite[Proposition 20.23]{bauschke:combettes}. 
Finally, strict monotonicity of $\Gamma^i = A \, \diag\left( (\partial f^i)^{-1}, (\partial f^0)^{-1} \right) \circ B$ follows from the next result, which is a variant of \cite[Proposition 28.2]{bauschke:combettes}.
\smallskip
\begin{lemma} \label{lem:monotonicity-preservation}
Let $\mathcal{M}: \R^m \rightarrow \R^m$ be a (strictly) monotone mapping, and $A, B \in\R^{m \times m}$. If $A$ is invertible and 
$A^{-\top} B \in \mathbb{S}_{\succ 0}^{m}$, then the mapping $A \, \mathcal{M}\left( B \, \cdot \right)$ is (strictly) monotone in $\mathcal{H}_{A^{-\top} B}$.
{\hfill $\square$}
\end{lemma}
\smallskip
\begin{proof}
{Since} $\mc{M}$ is (strictly) monotone, for all $x\neq y \in \R^m$, {we have}:
\begin{align*}
0 \leq (<) \, & \left( \mathcal{M}( B x ) - \mathcal{M}( B y ) \right)^\top \left( B x - B y \right) \\
& = \left( \mathcal{M}( B x ) - \mathcal{M}( B y ) \right)^\top 
A^{\top} A^{-\top} B \left( x-y \right)  \\
& = \left( A \, \mathcal{M}( B x ) - A \, \mathcal{M}( B y ) \right)^\top 
A^{-\top} B \left( x-y \right).
\end{align*}
\end{proof}

\subsection{Dynamic control as fixed point iteration: Proof of Theorem \ref{th:convergence-rate} (Global logarithmic convergence rate)}
\label{sec:proof-linear-rate}

The iteration in \eqref{eq:forward-backward} {can be written} as {the} fixed point iteration 
$$
\left[
\begin{matrix}
\sigma_{(t+1)} \\
\lambda_{(t+1)}
\end{matrix}
\right] = 
\left( 1-\alpha_t\right) \left[
\begin{matrix}
\sigma_{(t)} \\
\lambda_{(t)}
\end{matrix}
\right]
+ \alpha_t \, 
\mc{T}\left( 
\left[
\begin{matrix}
\sigma_{(t)} \\
\lambda_{(t)}
\end{matrix}
\right]
\right)
$$
where the mapping $\mathcal{T}: \R^{2n} \rightarrow \R^{2n}$ is defined as
\begin{equation} \label{eq:T}
\mc{T}(\cdot) := \left( I + \epsilon \, M \right)^{-1}\left( \textup{Id} - \epsilon \, \Gamma \right)(\cdot).
\end{equation}

In fact, 
a vector $\bar{z} \in \R^{2n}$ is a fixed point of 
$\mc{T}$ if and only if, for all $\alpha > 0$, it is fixed point of the mapping $(1-\alpha) \textup{Id} + \alpha \mc{T}$, and if and only if it is a zero of $M + \Gamma$ \cite[Lemma $1$]{grammatico:16netgcoop}: if $\bar{z} = \mathcal{T}(\bar z) = \left( I + \epsilon \, M \right)^{-1}\left( \textup{Id} - \epsilon \, \Gamma \right)(\bar{z})$, 
then $\left( I + \epsilon \, M \right) \bar{z} = \bar{z} - \epsilon \, \Gamma(\bar{z})$, which is equivalent to $M \bar{z} + \Gamma(\bar{z}) = \bs{0}$.

To establish the convergence rate of the iteration in \eqref{eq:forward-backward} with $\alpha_t = \bar{\alpha} \in (0,1]$ for all $t \in \N$, we show that the mappings $\mathcal{T}$ in \eqref{eq:T} 
and 
\begin{equation} \label{eq:K}
\mathcal{K}(\cdot) := (1-\bar\alpha) \textup{Id}(\cdot) + \bar{\alpha} \mathcal{T}(\cdot)
\end{equation}
are averaged operators.
\smallskip
\begin{lemma} \label{lem:T-averaged}
The mapping $\mathcal{T}$ in \eqref{eq:T} is $\frac{2}{3}$-averaged in $\mathcal{H}_{P}$, and the mapping $\mc{K}$ in \eqref{eq:K} is $(1\!-\!\frac{\bar{\alpha}}{3})$-averaged in $\mathcal{H}_{P}$, with $P$ as in \eqref{eq:P-C-K}.
{\hfill $\square$}
\end{lemma}
\smallskip
\begin{proof}
It follows from the proof of Theorem \ref{th:forward-backward} that $M$ in \eqref{eq:M-in-Theta} is monotone in $\mc{H}_{P}$, thus $(I + \epsilon M)^{-1}$ is firmly nonexpansive \cite[Corollary 23.10 (i)]{bauschke:combettes}.
According to Theorem \ref{th:monotone}, $\epsilon \, \Gamma$ is firmly nonexpansive, {hence also} the mapping $\textup{Id} - \epsilon \, \Gamma$ is firmly nonexpansive \cite[Proposition 4.2]{bauschke:combettes}. 
Therefore, $\mc{T}$ is the composition of two firmly nonexpansive mappings, or equivalently the composition of two $\frac{1}{2}$-averaged operators. In particular, $\mc{T}$ is $\frac{2}{3}$-averaged in $\mc{H}_P$ \cite[Proposition 2.4]{combettes:yamada:15}. 
Finally, $\mc{K}$ is the convex combination of two averaged operators, hence it is averaged with parameter $1-\bar{\alpha} + \bar{\alpha}\frac{2}{3} = 1- \frac{\bar{\alpha}}{3}$ \cite[Proposition 2.2]{combettes:yamada:15}.
\end{proof}
\smallskip

We can now prove Theorem \ref{th:convergence-rate}.
\begin{proof} 
By the definition of averaged operator, we have that
$ \left\| \mathcal{K}(x) - \mc{K}(y) \right\|_P^2 \leq \left\| x-y \right\|_P^2 - \frac{1-\beta}{\beta}\left\| \mc{K}(x) - \mc{K}(y) -(x-y)  \right\|_P^2 $ for all $x, y \in \R^{2n}$, where $\beta = 2/3$ in view of Lemma \ref{lem:T-averaged}. Let us take $x = z_{(\tau)} := \left[ \sigma_{(\tau)}; \lambda_{(\tau)} \right]$, and $y = \bar{z} = \mathcal{K}(\bar{z})$. By substituting, we obtain
$ \left\| \mathcal{K}\left(z_{(\tau)}\right) - \bar{z} \right\|_P^2 \leq \left\| z_{(\tau)} - \bar{z} \right\|_P^2 - \frac{1-\beta}{\beta} 
\left\| \mc{K}\left(z_{(\tau)}\right) - z_{(\tau)} \right\|_P^2$, and equivalently 
$\left\| \mc{K}\left(z_{(\tau)}\right) - z_{(\tau)} \right\|_P^2
\leq  \left( \frac{\beta}{1-\beta}  \right) \left\| z_{(\tau)} - \bar{z} \right\|_P^2 - \left\| \mathcal{K}\left(z_{(\tau)}\right) - \bar{z} \right\|_P^2$. In particular, note that 
$\left\| \mc{K}\left(z_{(\tau)}\right) - \bar{z} \right\|_P^2 \leq \left\| z_{(\tau)} - \bar{z} \right\|_P^2$.
Now, we sum up over $\tau \in \N[0,t]$ and derive 
\begin{align*}
 & (t+1) \left\| z_{(t+1)} - z_{(t)} \right\|_P^2 \\
 &\leq 
\textstyle \sum_{\tau = 0}^t \left\| z_{(\tau)} - z_{(\tau)} \right\|^2 \\
 &\leq \textstyle \frac{\beta}{1-\beta} \sum_{\tau = 0}^t \left\| z_{(\tau)} - \bar{z} \right\|_P^2 - \left\| \mathcal{K}\left(z_{(\tau)}\right) - \bar{z} \right\|_P^2 \\
  & \leq \textstyle \frac{\beta}{1-\beta}
  \left\| z_{(0)} - \bar{z} \right\|_P^2.
\end{align*} 
Since $\beta = 1-\frac{\bar{\alpha}}{3}$, we have that $\left\| z_{(t+1)} - z_{(t)} \right\|_P^2 \leq  
\frac{1}{t+1} \left( \frac{3}{\bar{\alpha}}-1\right) \left\| z_{(0)} - \bar{z} \right\|_P^2$, which completes the proof.
\end{proof}
\smallskip

Finally, we note that the mapping $\mathcal{T}$ is nonexpansive according to Lemma \ref{lem:T-averaged}, hence 
several fixed point iterations have global convergence guarantee. 
The design choice $\alpha_t := \bar{\alpha} \in (0,1)$ for all $t \in \N$ is known as Krasnoselskij iteration \cite[Equation (18)]{grammatico:parise:colombino:lygeros:16}. Among other iterations, we mention the Mann iteration \cite[Equation (20)]{grammatico:parise:colombino:lygeros:16}, which corresponds to choosing the sequence 
$(\alpha_t)_{t\in \N}$ such that $\alpha_t \in [0,1]$ for all $t \in \N$, $\lim_{ t \rightarrow \infty} \alpha_t = \infty$, and  $\sum_{t=0}^{\infty} \alpha_t = \infty$, e.g. $\alpha_t := 1/(t+1)$ for all $t \in \N$.

\section{Discussion} \label{sec:discussion}

\subsection{Features of the dynamic control scheme}

One computational feature of the iteration in \eqref{eq:forward-backward} is that it only requires one-to-all coordination between a central computer and the decentralized, hence \textit{parallelizable}, optimal responses $( x^{i \star} )_{i=1}^{N}$ in \eqref{eq:optimal-responses} of the agents, as summarized in Algorithm \ref{alg:aggregative-control}. Each decentralized computation consists in solving a finite dimensional strongly convex optimization problem, for which efficient algorithms are available.
Note that at each iteration $t$ {only one vector in $\R^n$ is broadcast}, independently on the population size $N$, which can be arbitrarily large. Also note that the coordinator needs access to the aggregate information $\mathcal{A}$ only, not necessarily to the entire set of optimal responses.

The main distinctive feature of the proposed {semi-decentralized} architecture is that the central coordinator can decide on the step sizes $\left( \alpha_t \right)_{t=0}^{\infty}$ and $\epsilon$, on the gain $K$, and on the stopping criterion for the iteration in \eqref{eq:forward-backward}.
Therefore, the agents can simply behave as \textit{fully noncooperative}, in the sense that in addition to being self-interested, they do not have to exchange information with each other, nor to agree on the step sizes associated with the control signals, nor on the stopping criterion. Each of these agreement points would be in fact exposed to malicious agent behavior. Also, for the convergence of the dynamic {iterations}, neither the central coordinator nor the agents need to know the population size $N$.

In summary, Algorithm \ref{alg:aggregative-control} is such that:
\begin{itemize}
\item[$\bullet$] the central coordinator keeps the information on the chosen incentive mechanism and on the global coupling constraint private;
\item[$\circ$] the noncooperative agents keep information on their cost functions and local constraints private.
\end{itemize}

\subsection{Generalized Nash aggregative games}
The game setup in \eqref{eq:shared-constraint}--\eqref{eq:best-response-price} can be related to the generalized Nash equilibrium problem ({\small GNEP}) with best responses
\begin{align} \label{eq:GNEP}
\begin{split}
x^i \in \underset{y^i \in \R^n}{\argmin} & \
\textstyle \varphi^i\left( y^i, \, \bs{x}^{-i}\right)
\\
\textup{s.t.} & \ y^i \in \mathcal{X}^i \cap \mathcal{S}^i\left( \bs{x}^{-i}\right),
\end{split}
\end{align}
for all $i \in \N[1,N]$, where
$$ \textstyle \varphi^i\left( y^i, \, \bs{x}^{-i}\right) := 
f^i\left( y^i\right)  +  
\left( C \, \tfrac{1}{N} \left( y^i + \sum_{j \neq i}^N x^j  \right) \right)^\top \! y^i, $$
and the shared constraint set for agent $i$ reads as
$$\textstyle \mathcal{S}^i\left( \bs{x}^{-i} \right) := 
\left\{ y^i \in \R^n \mid \tfrac{1}{N} (y^i + \sum_{j \neq i}^{N} x^j) \in \mathcal{S}\right\}.$$
Due to the aggregative structure of both the cost functions and the shared constraint, let us label such a game as generalized Nash aggregative game.

Since there is one unique shared constraint that is convex, the game is called \textit{jointly convex} \cite[\S 3.2]{facchinei:kanzow:07}. 
Several methods are available for solving a jointly convex {\small GNEP}, e.g. the decomposition approach outlined in \cite[Part II, \S 3, p. 166]{cominetti:facchinei:lasserre} and summarized next. First, we shall assume that the shared constraint
$\textstyle \tfrac{1}{N} \sum_{i=1}^{N} x^i \in \mathcal{S} $ can be written as 
$$\textstyle g(\bs{x}) := \tfrac{1}{N}\sum_{i=1}^N g^i(x^i) \leq 0$$ 
for some convex, differentiable functions $g$ and $\{g^i\}_{i=1}^{N}$. Then, we introduce an additional agent that controls the dual variable $\lambda \in \R^m$ associated with the shared constraint, and let each agent minimize its own Lagrangian function.
Specifically, we derive the following (non-generalized) NEP among $N+1$ agents with no coupling constraint \cite[Part II, Equation (3.2)]{cominetti:facchinei:lasserre}:
\begin{align} \label{eq:GNEP-2}
\begin{split}
x^i \in \underset{y^i \in \mathcal{X}^i}{\argmin} & \
 \varphi^i\left( y^i, \, \bs{x}^{-i}\right) + \tfrac{1}{N} \, \lambda^\top \,  g^i\left( y^i\right) \\
\lambda \in \underset{\zeta \in \R_{\geq 0}^m}{\argmin} & \ - \! \zeta^\top g(\bs{x}).
\end{split}
\end{align}
It then follows from \cite[Theorem 8]{facchinei:kanzow:07} that, under basic regularity assumptions, an NE $\left( \bs{x}, \lambda \right)$ for the game without coupling constraint \eqref{eq:GNEP-2} is such that $\bs{x}$ is a {\small GNE} for the game with coupling constraint in \eqref{eq:GNEP}. 
We note that in \eqref{eq:GNEP-2} the best response for $x^i$ is similar to that in \eqref{eq:best-response-price}, while the best response for $\lambda$ has no clear counterpart in \eqref{eq:forward-backward}.

For the computation of an NE for \eqref{eq:GNEP-2} several distributed algorithms are available, see e.g. \cite[Part II, \S 2]{cominetti:facchinei:lasserre}, or follow from solution algorithms for monotone variational inequalities \cite[\S 12]{facchinei:pang}. We note that all such computational algorithms require differentiable cost and coupling constraint functions and that the so-called game (or pseudo gradient) mapping
$$\mc{F} \left(\bs{x}, \lambda \right) :=  
\left[
\begin{matrix}
\partial_{x^1} \varphi^1\left( x^1, \, \bs{x}^{-1} \right) + \tfrac{1}{N} \,  \partial g^1 \!\left( x^1\right) \lambda \\
\vdots\\
\partial_{x^N} \varphi^N\left( x^N, \, \bs{x}^{-N} \right) + \tfrac{1}{N} \,  \partial g^N \!\left( x^N\right) \lambda \smallskip \\
-g(\bs{x})
\end{matrix}
\right]
$$
for \eqref{eq:GNEP-2} is (strongly) monotone in some Hilbert space. However, monotonicity of the game mapping does not hold in general \cite[\S 5.2]{facchinei:kanzow:07}, not even with convex problem data. Even in the aggregative game setup in \eqref{eq:shared-constraint}--\eqref{eq:best-response-price}, or in  \eqref{eq:GNEP}, Standing Assumptions \ref{ass:standing}--\ref{stand:J} do not imply that the game mapping is monotone either.

We finally mention that, since $C \in \bbS^n$, it can be shown that the {\small GNEP} in \eqref{eq:GNEP} is a generalized potential game (GPG) \cite[Definition 2.1]{facchinei:piccialli:sciandrone:11}. 
Thus, regularized Gauss--Siedel algorithms are applicable for computing a GNE, e.g. \cite[Algorithms 2, 3]{facchinei:piccialli:sciandrone:11}, as monotonicity of the game mapping is not required.

\subsection{Generalized projected dynamical systems}

The problem to control the decentralized optimal responses of the agents to an equilibrium in Section \ref{sec:game} can be {interpreted} as the design of a dynamic output-feedback control law 
{$u(t) = \left[ C, \, K\right] \kappa(t, y(t))$} 
for the following discrete-time system:
\begin{subequations}
\begin{align}
x^i(t+1) & = \underset{\xi \in \mc{X}^i}{\argmin} \ 
f^i\left( \xi \right) + \xi^\top u(t), \ \ \forall i \in \N[1,N], 
\label{eq:argmin-dynamics}
\\
y(t) & = \textstyle  \frac{1}{N} \sum_{i=1}^{N} x^i(t).
\end{align}
\end{subequations}
In particular, the control objective is to drive the output $y$ to some $\bar{y} \in \mc{S}$ and the states to an equilibrium $\bar{\bs{x}} \in (\mc{X}^1 \times \ldots \times \mc{X}^N)$ such that, for all $i \in \N[1,N]$, 
$\bar{x}^i = {\arg\min}_{\xi \in \mc{X}^i} \, f^i\left( \xi\right) + \bar{u}^\top \xi$, for some $\bar{u} \in \R^n$.

Whenever the functions {$\{ f^i \}_{i=1}^{N}$} are strongly convex quadratic, the $\argmin$ dynamics in \eqref{eq:argmin-dynamics} read as \textit{projection dynamics}. 
Precisely, if $f^i(\xi) := \frac{1}{2}\xi^\top Q^i \xi + {c^i}^\top \xi$ for some $Q^i \succ 0$ and $c^i \in \R^n$,
then we have 
$x^{i}(t+1) = \proj_{\mc{X}^i}^{Q^i}( -{Q^i}^{-1}( c^i + u(t)) )$, 
where the projection operator is defined as $\proj_{\mc{X}^i}^{Q^i}(\cdot) := \arg\min_{\xi \in \mc{X}^i} \left\| \xi - \cdot \right\|_{Q^i}$.
Consequently, \eqref{eq:argmin-dynamics} is a discrete-time projected dynamical system \cite{dupuis:nagurney:93, nagurney:zhang} where we wish to close the feedback loop with some 
$u(t) = \left[ C, \, K\right] \kappa\left(t, \frac{1}{N} \sum_{i=1}^{N} x^i(t)\right)$.
For general strongly convex functions as in Standing Assumption \ref{stand:J}, it follows from Lemma \ref{lem:optimizer} that in \eqref{eq:argmin-dynamics} we have 
$x^i(t+1) = (\partial f^i)^{-1} \!\left( -u(t) \right)$, 
where $(\partial f^i)^{-1}$ is strictly monotone and Lipschitz continuous according to Lemma \ref{lem:strong-convexity-monotonicity}. Since the stability analysis for projected dynamical systems builds upon the fact that the projection onto a convex set is a firmly nonexpansive (hence strictly monotone and $1$-Lipschitz continuous) mapping, our analysis based on monotone operator theory is applicable to generalized projected dynamical systems with strongly convex \textit{non-quadratic} functions.

\section{Applications} \label{sec:applicability}
The {considered} aggregative game setup with coupling constraints is applicable to the dynamic management of noncooperative agents coupled in aggregative form. Applications include network congestion control \cite{barrera:garcia:15}, demand response in competitive  markets \cite{li:chen:dahleh:15} and demand side management for prosumers in the smart grid, e.g. residential loads with coupling constraints \cite{deng:yang:chen:asr:chow:14, deng:xiao:lu:chen:15}, and smart homes with shared renewable energy sources \cite{carli:dotoli:15}. The common feature of all these setups is in fact the presence of a population of competitive agents with convex cost functions, convex local and shared constraints, coupled together in aggregative form as in \eqref{eq:shared-constraint}--\eqref{eq:best-response-price}. In the next subsections, we focus on two such applications.

\subsection{Network congestion control with capacity constraints}
We consider the problem faced by a network manager to control the flow demands of a large set of noncooperative users by dynamically pricing the network capacity \cite{barrera:garcia:15}.

\subsubsection*{General problem setup}
The network is characterized by a set of edges $\mc{E} := \N[1,n]$, with capacity $\overline{c} := [ \overline{c}_1; \, \ldots; \, \overline{c}_n ] \in \R_{>0}^{n}$. Each user $i \in \N[1,N]$ aims at selecting its flow profile $x^i \in \R^n$ that minimizes its disutility function, that is coupled in aggregative form to the flow profiles of all other agents. Specifically, each user $i$ aims at minimizing the cost function 
\begin{equation}\label{eq:Ji-congestion}
f^i\left(x^i\right) + \textstyle \left( {c}\left( \frac{1}{N} \sum_{j=1}^{N} x^j \right) + \lambda \right)^\top x^i,
\end{equation}
where \cite[Sections III-IV]{barrera:garcia:15} the function $f^i$, that is continuous and convex, represents the intrinsic disutility, the function $c(\cdot) = [c_1(\cdot); \, \ldots; \, c_n(\cdot)]^\top$ represents the flow-unit delay cost (disutility) experienced by the users over the edges, $\lambda \in \R^n$ is the congestion price, that is, the penalty vector associated with the coupling network-capacity constraint
\begin{equation}\label{eq:shared-congestion}
\bs{0} \leq \textstyle \frac{1}{N}\sum_{i=1}^{N}{x^i} \leq \overline{c}.
\end{equation} 
Each user $i$ also has a local constraint set $\mc{X}^i$ that  represents its individual routing policy per flow unit.

\subsubsection*{Illustrative scenario with fixed routing policy}
In the following, we simulate the scenario illustrated in \cite[Section IV.B]{barrera:garcia:15}, where fixed routing policies are considered. Namely, each agent $i$ has a routing policy $x^i = a^i \xi^i$, for some scalar $\xi^i \geq 0$, and some fixed vector $a^i \in \R_{\geq 0}^n$ such that $\bs{1}^\top a^i = 1$. 
As in \cite[Section IV.B]{barrera:garcia:15}, for each user $i$, we use the convex, non-quadratic, intrinsic disutility function $f^i( \xi^i ) := - 20 \, \textup{ln}( 1 + \xi^i )$.

\subsubsection*{Affine approximation of the delay cost}
Next we derive an affine approximation for the delay mapping $c(\cdot)$, based on the function
$c_e( \cdot ) = 1/(\beta_e - \cdot)$ from \cite{barrera:garcia:15}, for each edge $e \in \N[1,N]$.
We use the first order Taylor approximation around the origin, hence in \eqref{eq:Ji-congestion} we consider the delay cost mapping
\begin{equation*}
c \left( \cdot \right) := \textstyle
\left[ \diag \left(\left( {1}/{\beta_e^2} \right)_{e=1}^{n} \right)
\right] \left( \cdot \right) + \vect\left(  \left( {1}/{\beta_e} \right)_{e=1}^{n} \right).
\end{equation*}
Note that the matrix $\diag \left(\left( {1}/{\beta_e^2} \right)_{e=1}^{n} \right) \succ 0$ corresponds to the matrix $C$ in the cost-function structure in \eqref{eq:J}, hence any matrix gain $K \succ 0$ satisfies the design choice \ref{des:K}. 

\subsubsection*{Numerical simulations}
\begin{figure}
\begin{center}
\includegraphics[width=1\columnwidth]{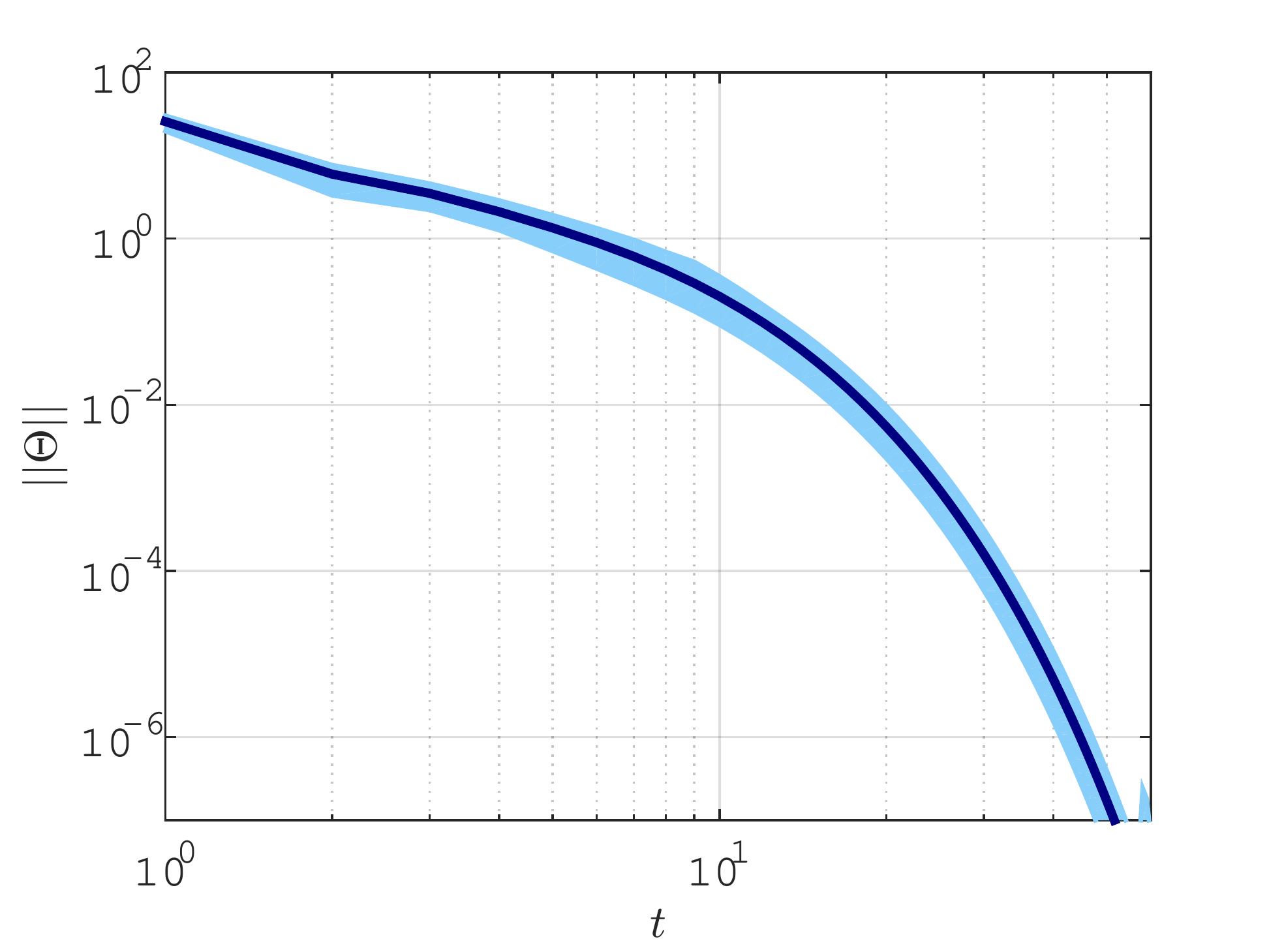}
\end{center}
\caption{
$\left\| \Theta\left( [ \sigma_{(t)}; \lambda_{(t)} ]\right) \right\|$ in $10^3$ experiments, each with $N = 10^{4}$ agents, as a function of the iteration number $t$.
The shaded area represents all the experiments; the solid line is the average value.
}
\label{fig:all_theta}
\end{figure}

\begin{figure}
\begin{center}
\includegraphics[width=1\columnwidth]{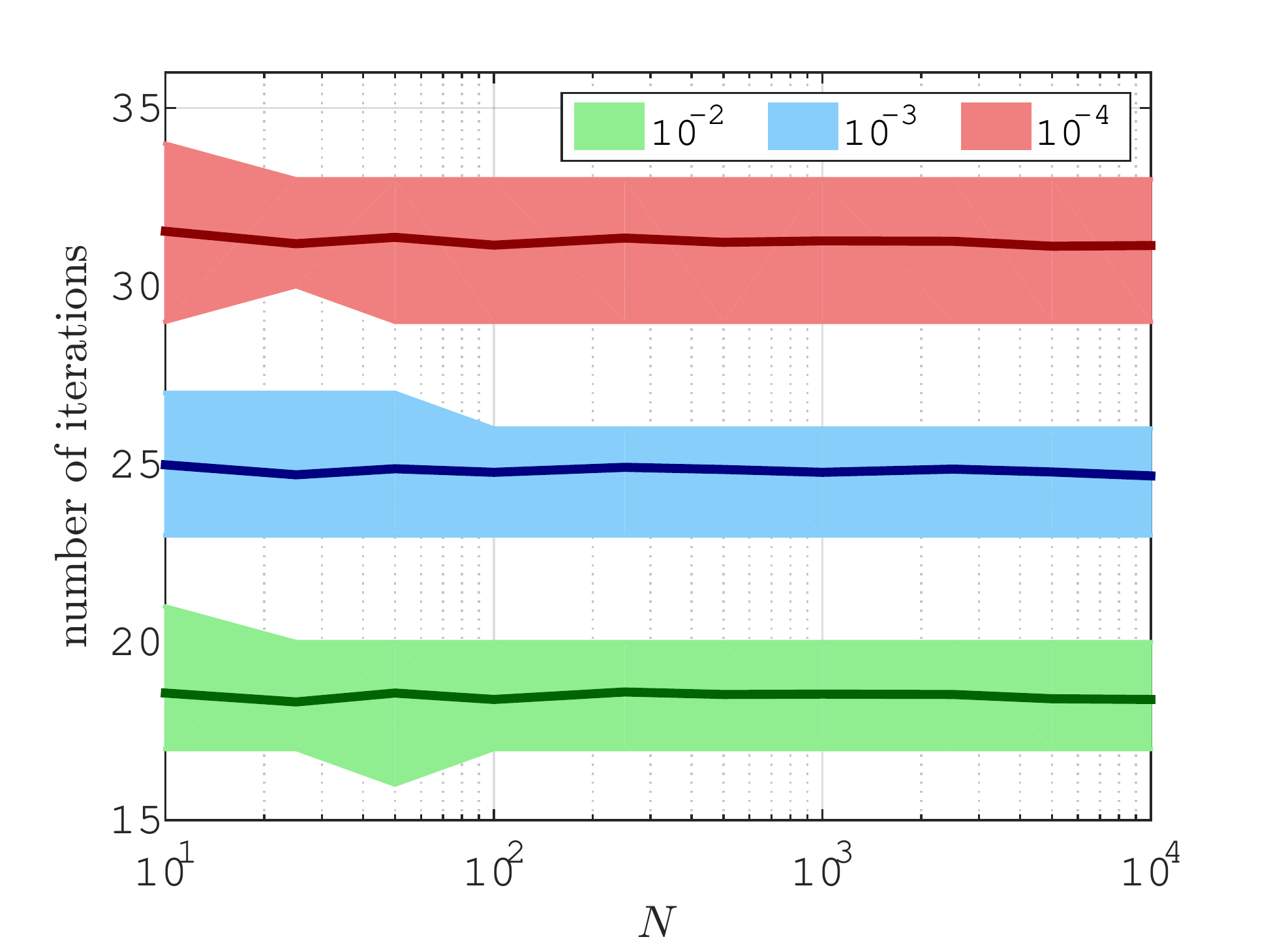}
\end{center}
\caption{
Number of iterations needed for 
$\left\| \Theta\left( [ \sigma_{(t)}; \lambda_{t} ] \right)\right\|$  to be less than $10^{-2}$ (green), $10^{-3}$ (blue), 
$10^{-4}$ (red), as a function of the population size $N$. 
The shaded areas represent the number of iterations for the whole set of experiments; the solid lines represent the average number of iterations.
}
\label{fig:full_convergence}
\end{figure}

We use some numerical parameters from \cite[Tables I, II]{barrera:garcia:15}, namely $n=5$, $\beta_1 = \ldots = \beta_n = 20$, together with capacities $\overline{c}_2 = \ldots = \overline{c}_{4} = 4$, $\overline{c}_{1} = \overline{c}_{5} = 2$, and stochastic vectors $a_1, \ldots, a_N$ sampled with uniform distribution. The local constraints are set as $\bs{0} \leq x^i \leq 10 \cdot \bs{1}$ for all $i \in \N[1,N]$.
We tune the gains of the dynamic control law $\kappa$ in \eqref{eq:forward-backward} as $\alpha=1$, $K= I$, and choose 
$\epsilon$ that satisfy the design choice \ref{des:epsilon}. We take as initial condition $\lambda_{(0)} = 0$ and  random $\sigma_{(0)}$, uniformly distributed within the coupling constraint set. Finally, as convergence criteria, we consider that $\left\| \Theta\left( [\sigma_{(t)}; \lambda_{(t)}] \right) \right\|$ must be less than certain tolerance values.
We run several numerical experiments, each with randomly selected routing policies $\{a^i\}_{i=1}^{N}$ and initial condition $\sigma_{(0)}$. Figure \ref{fig:all_theta} shows the convergence scenario for $10^3$ experiments, with $N = 10^4$ agents. Figure \ref{fig:full_convergence} shows the convergence scenarios parametric on the population size, where $10^2$ experiments are run for each value of $N$. We conclude that the population size does not affect the convergence speed.

\subsection{Charging coordination for plug-in electric vehicles with transmission line constraints}
We consider the problem to control the charging schedule of a population of plug-in electric vehicles subject to transmission line constraints \cite{ghavami:kar:bhattacharya:gupta:13, tushar:saad:poor:smith:12}.

\subsubsection*{Problem setup}
Each user $i$ aims at charging its vehicle with energy injections $[ x_1^i, x_2^i, \ldots, x_n^i] =: x^i$,
while minimizing its individual disutility, subject to individual and shared charging constraints, over a charging horizon that is here discretized into $n=14$ time intervals. 
The nominal values of the numerical parameters defining the cost functions and the charging constraints are taken from \cite{ma:zou:liu:15}, and then are randomized as in \cite{grammatico:16cdc-pev} to emulate the population variability.
For each PEV agent $i \in \N[1,N]$, we consider the quadratic cost function $J^i(x^i, \sigma) = q^i \, {x^i}^\top x^i + {c^i}^\top x^i + (a \,  \sigma + b \bs{1}_n)^\top x^i$ that represents the battery degradation cost $q^i \, {x^i}^\top x^i + {c^i}^\top x^i$ \cite[Section II.C]{ma:zou:liu:15}, 
plus the electricity pricing $(a \, (\sigma+d) + b \bs{1}_n)^\top x^i$, where $a>0$ represents the inverse of the price elasticity of demand, $b>0$ represents the baseline price, and the vector $d \in \R^n$ the normalized average inflexible demand.

\subsubsection*{Numerical parameters}
With uniform distribution, we sample $q^i \sim \{0.004\} + [-0.002, \, 0.002]$ and $c^i \sim \{ 0.075 \} + [-0.02, \, 0.02]$. Further, we consider the normalized charging constraints $x^i \in \mathcal{X}^i := [\bs{0}_n, \overline{x}^i] \cap \left\{ y \in \R^n \mid \bs{1}_n^\top y = \gamma^i \right\}$, where with uniform distribution we sample $\gamma^i \sim \{ 0.8 \} + [-0.2, \, 0.2]$, and the vector $\overline{x}^i \in \R^n$ is such that, for all $j \in \N[1,n]$, $\overline{x}_j^i \sim \{ 0, 0.25 \}$, with $\overline{x}_j^i = 0$ (that is, no charging at the time interval $j$) with probability $20\%$. In addition, for $20\%$ of the overall population,  we consider the vehicle-to-grid option, namely by substituting the lower and upper bounds $[\bs{0}_n, \overline{x}^i]$ with $[-\frac{1}{2} \overline{x}^i, \overline{x}^i]$.
Next, we scale the parameters in \cite[Section IV]{ma:zou:liu:15} with respect to the population size, and derive the parameters $a = 0.038$, $b = 0.06$ and $d$ empirically from \cite[Figure 1]{ma:callaway:hiskens:13}, \cite[Figure 1]{ma:zou:liu:15}.

\subsubsection*{Coupling constraints}
We extend the setup in \cite{grammatico:16cdc-pev} with time-varying transmission line constraints, i.e., $\bs{0}_n \leq \frac{1}{N} \sum_{i=1}^{N} x^i \leq \overline{c}$ as in \eqref{eq:shared-congestion}
\cite[Equation (13)]{ghavami:kar:bhattacharya:gupta:13},
\cite[Equation (1)]{tushar:saad:poor:smith:12}, for some vector $\overline{c} \in \R_{\geq 0}^n$. We illustrate the proposed algorithm with capacities $\overline{c}_j = 0.04$ if $j \in \{1, 2, 12, 13, 14\}$, $0.1$ otherwise, {to represent more restrictive} charging limitations during the day time.

\subsubsection*{Numerical simulations}
\begin{figure}
\centering
\includegraphics[width=1\columnwidth]{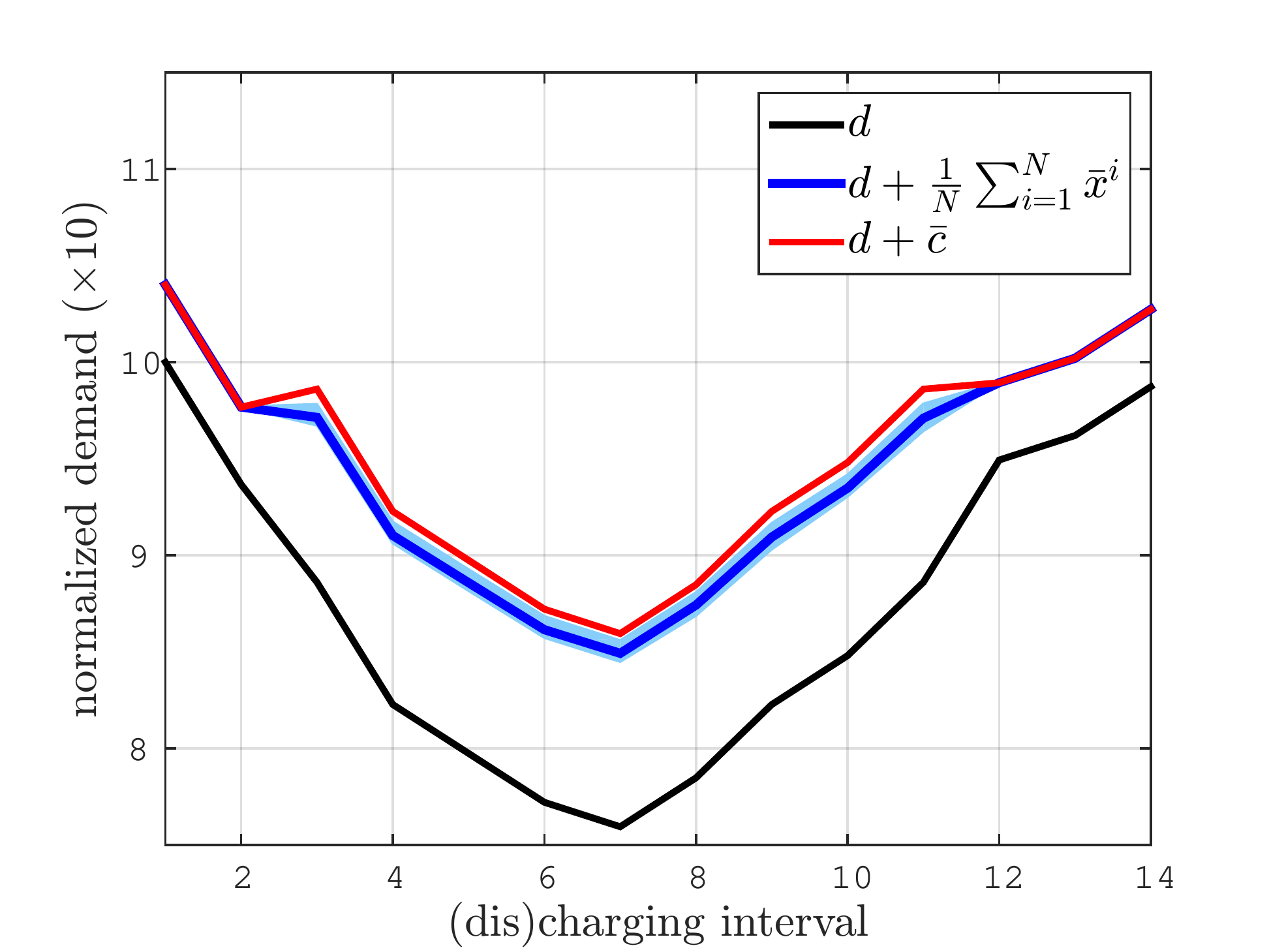}
\caption{Sum between the normalized average inflexible demand $d$ and the average among the charging strategies $\frac{1}{N} \sum_{i=1}^{N} x^{i \star}( \bar{\sigma}, \bar{\lambda} ) $ at the equilibrium. The shaded area represents the union over all the experiments.
}
\label{fig:valley}
\end{figure}

\begin{figure}
\begin{center}
\includegraphics[width=1\columnwidth]{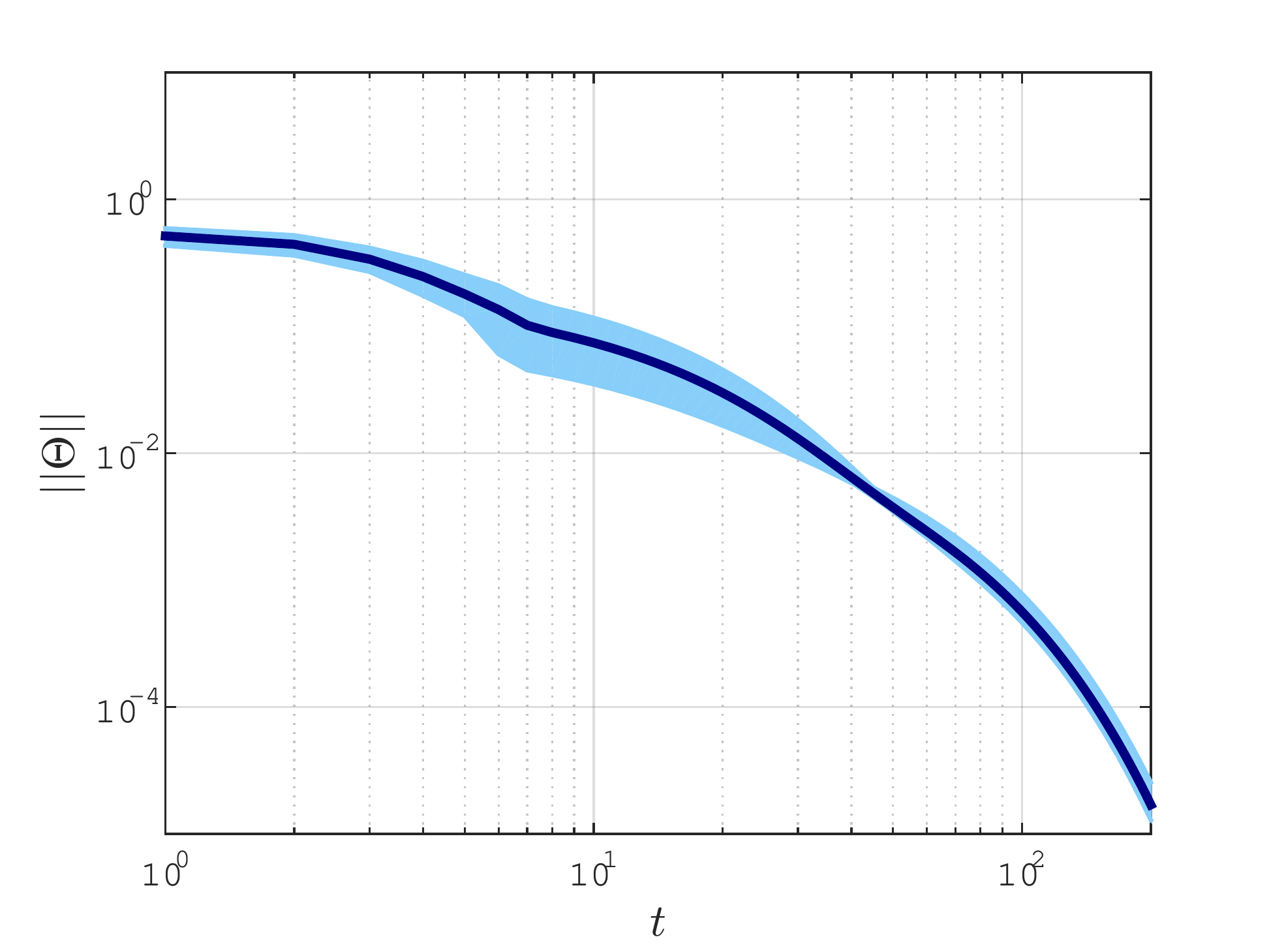}
\end{center}
\caption{
$\left\| \Theta\left( [ \sigma_{(t)}; \lambda_{(t)} ]\right) \right\|$ in $10^3$ experiments, each with $N = 10^{4}$ agents, as a function of the iteration number $t$.
The shaded area represents the union over all the experiments; the blue solid line is the average value.
}
\label{fig:all_theta_pev}
\end{figure}

\begin{figure}
\begin{center}
\includegraphics[width=1\columnwidth]{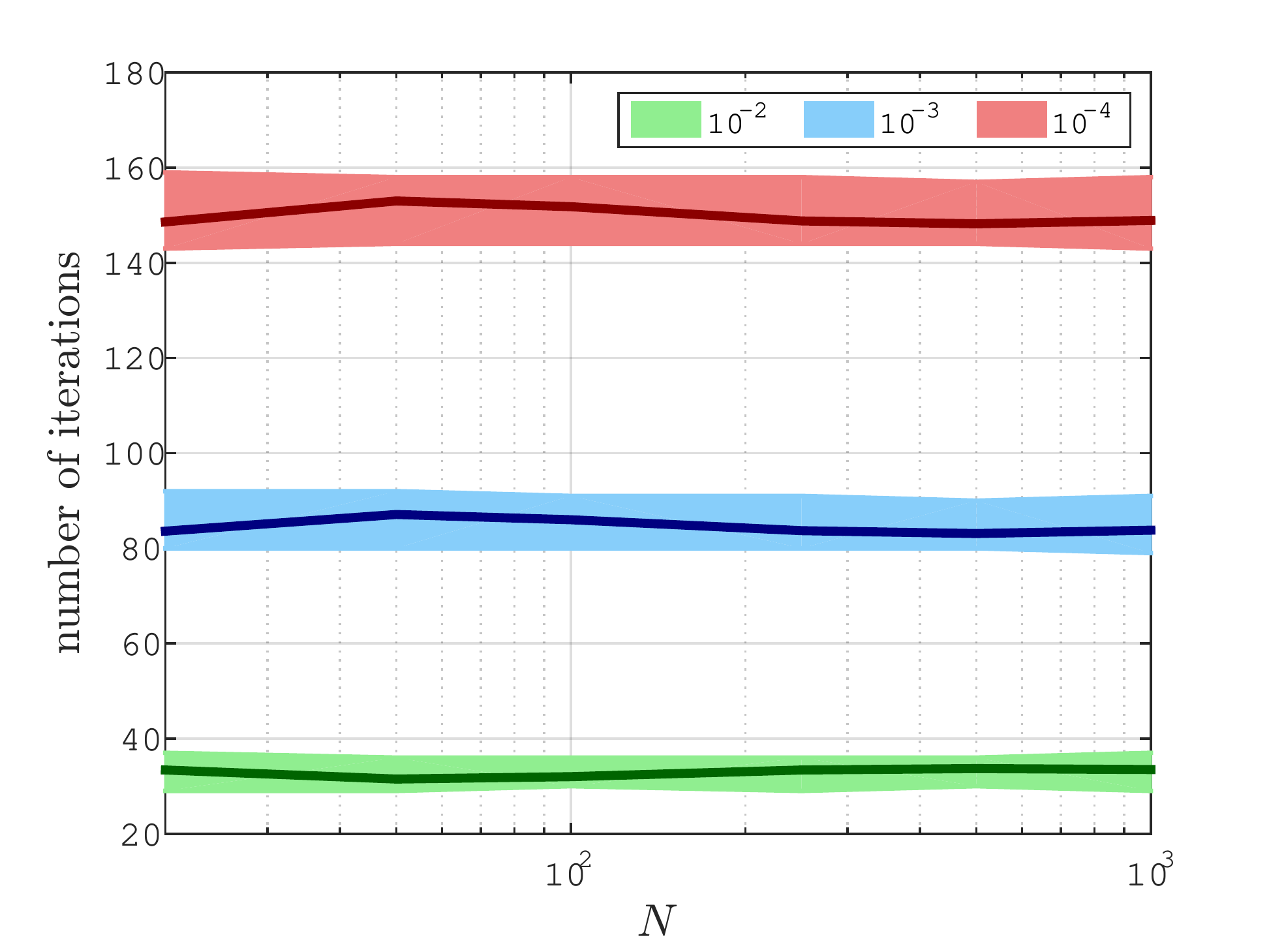}
\end{center}
\caption{
Number of iterations needed for 
$\left\| \Theta\left( [ \sigma_{(t)}; \lambda_{t} ] \right)\right\|$  to be less than $10^{-2}$ (green), $10^{-3}$ (blue), 
$10^{-4}$ (red), as a function of the population size $N$. 
The shaded areas represent the number of iterations for the whole set of the experiments; the solid lines represent the average number of iterations.
}
\label{fig:full_convergence_pev}
\end{figure}

We tune the gains of the dynamic control law in \eqref{eq:forward-backward} as $\alpha=1$, $K=0.05 I$, and then choose $\epsilon$ according to the design choice \ref{des:epsilon}. 
We take as initial condition $\lambda_{(0)} = 0$ and random $\sigma_{(0)}$, uniformly distributed within the shared constraint set. Finally, as convergence criteria we consider that $\left\| \Theta\left( [\sigma_{(t)} ; \lambda_{(t)}] \right) \right\|$ shall be less than certain tolerance values.
We run several numerical experiments, each with the mentioned randomly selected parameters and initial condition. 
Figure \ref{fig:valley} shows the sum between the normalized average inflexible demand and the average among the charging strategies at the equilibrium, i.e., the optimal responses in \eqref{eq:optimal-responses}, for $10^2$ experiments with $N = 10^4$ agents;
Figure \ref{fig:all_theta_pev} shows the convergence scenario for $10^3$ experiments. Figure \ref{fig:full_convergence_pev} shows the convergence scenarios parametric on the population size, where $10^2$ experiments are run for each value of $N$. Also for this application, we conclude that the population size does not affect the convergence speed of the proposed algorithm.

\section{Conclusion and outlook} \label{sec:conclusion}
We have addressed the problem to control a population of competitive agents, with convex cost functions coupled together via the average population state, convex local and coupling constraints, towards an aggregative equilibrium. Our technical results allow us to design a model-free dynamic control law with global convergence guarantee, with no assumption on the problem data, other than strong convexity and compactness. The numerical simulations show that the proposed algorithm achieves an aggregative equilibrium within a reasonable number of iterations independently on the population size.

This paper can be extended in several directions, including those outlined next.
\begin{itemize}
\item \textit{Stochasticity}: We have not considered agents with probabilistic constraints. Such an extension would be valuable for the analysis and control of decentralized optimal decisions under stochastic uncertainty. One possible approach is to solve a deterministic approximation of the stochastic game, and then conclude about the computed solution in probabilistic terms.
 
\item \textit{Dynamic games}: We have considered agents that update their decision vector based on the signals received from the coordinator, but not on their decision history. More generally, it would be valuable to analyze the dynamics of agents with memory and cumulative cost functions. 

\item \textit{Asynchronous updates}:  
A framework with agents that update their strategies 
asynchronously would be more general and practically relevant in networked applications.

\item \textit{Optimal parameter selection}:
We have presented feasible design choices for the control parameters. To maximize the convergence speed, the parameters shall be appropriately selected, e.g. optimally with respect to some convergence rate estimate.

\end{itemize}

The connection with multi-agent dynamics in cooperative optimization is currently an active area of research \cite{li:zhang:zhao:lian:kalsi:16,deori:margellos:prandini:17}.

\appendix

\subsection{Proof of Proposition \ref{prop:existence-ANESC}}
\label{app:existence-equilibrium}

For all $\bs{\sigma} := \left[ \sigma^1 ; \ldots ; \sigma^N \right] \in \R^{n N}$, we consider the optimization problem
\begin{equation} \label{eq:P-bold-sigma}
\mc{P}( \bs{\sigma} ): \ 
\left\{ 
\begin{array}{cl}
\displaystyle \min_{ y^0 , \, \bs{y} } & \tilde{J}^{0}(y^0) + \sum_{i=1}^{N} \tilde{J}^i\left( y^i, \bs{\sigma} \right) \\
\textup{s.t.} & {K} \left(N \, y^0 - \sum_{j=1}^{N} y^j \right) = 0
\end{array}
\right.
\end{equation}
where $y^0 \in \R^n$, $\bs{y} := [ y^1; \ldots; y^N ] \in \R^{nN}$, $\tilde{J}^0: \R^n \rightarrow \overline{\R}$ is a strongly convex function such that $\textup{dom}(\tilde{J}^0) = \mc{S}$, and for all $i \in \N[1,N]$, $\tilde{J}^i: \R^n \times \R^{n N} \rightarrow \overline{\R}$ is defined as 
$\tilde{J}^i\left( y^i, \bs{\sigma}\right) := J^i( y^i , \frac{1}{N} \sum_{j =1 }^N \sigma^j, 0 )$. 

Let $\left[ x^{\star 0}(\bs\sigma); \, \bs{x}^{\star}(\bs\sigma) \right] = \left[ x^{\star 0}(\bs\sigma); x^{\star 1}(\bs\sigma) ; \ldots ; \, x^{\star N}(\bs\sigma) \right]$ be the optimizer of $\mc{P}( \bs{\sigma} )$ in \eqref{eq:P-bold-sigma}. 
Since the mapping $\bs{x}^{\star}(\bs{\cdot})$ takes values in the compact set $\bs{\mc{X}} :=  \mathcal{X}^1 \times \ldots \times \mathcal{X}^N $ and is (Lipschitz) continuous by Lemma \ref{lem:optimizer}, it has at least one fixed point \cite[Theorem 4.1.5 (b)]{smart}.
In the remainder of the proof, let $\bar{\bs\sigma} = \bs{x}^{\star}(\bar{\bs\sigma}) \in \bs{\mc{X}}$ be a fixed point of $\bs{x}^{\star}$.

Since $\mc{P}(\bar{\bs\sigma})$ from \eqref{eq:P-bold-sigma} has a separable convex cost function and the linear coupling constraint $y^0 = \frac{1}{N} \sum_{i=1}^{N} y^i$, see \cite[Equation 2.1]{boyd:admm}, it can be solved via the dual decomposition method \cite[Section 2]{boyd:admm}.
Specifically, the Lagrangian function 
\begin{align*}
\mc{L}\left( \bs{y}, \lambda \right) := & \ \textstyle  \tilde{J}({y^0}) + \sum_{i=1}^{N} \{ \tilde{J}^i\left( y^i, \bar{\bs\sigma} \right) \} \\ 
& \ + \lambda^\top {K} \textstyle \left( -N y^0 + \sum_{j=1}^{N} y^j\right) \\
\, = & \ \tilde{J}({y^0}) - N \lambda^\top {K} y^0
+ \textstyle \sum_{i=1}^{N} \tilde{J}^i\left( y^i, \bar{\bs\sigma} \right) + \lambda^\top {K} y^i
\end{align*}
is separable, therefore, due to Standing Assumption \ref{ass:standing} \cite[Section 5.2.3]{boyd:vandenberghe}, the iteration \cite[Equations 2.4, 2.5]{boyd:admm}
\begin{subequations} \label{eq:dual-decomposition-bold}
\begin{align}
x_{(t+1)}^{0} := & \, \underset{ y \in \R^n }{\argmin} 
\, \tilde{J}^0\left( y\right) - N \lambda_{(t)}^\top {K} y, 
\label{eq:dual-decomposition-bold-0} \\
x_{(t+1)}^{i} := & \, \underset{ y \in \R^n }{\argmin} 
\, \tilde{J}^i\left( y, \bar{\bs\sigma} \right) + \lambda_{(t)}^\top {K} y, \quad \forall i \in \N[1,N], 
\label{eq:dual-decomposition-bold-1} \\
\lambda_{(t+1)} := & \textstyle \, \lambda_{(t)} + \epsilon_t \, {K} \left( - N x_{(t+1)}^0 + \sum_{i=1}^N x_{(t+1)}^i \right)
\label{eq:dual-decomposition-bold-2}
\end{align}
\end{subequations}
converges to $\left( x^{ \star 0 }(\bar{\bs\sigma}), \bs{x}^{\star}(\bar{\bs\sigma}), \lambda_{\bar{\bs\sigma}}^{\star} \right)$, 
for an opportune choice of the sequence $\left( \epsilon_t \right)_{t=0}^{\infty}$, where $\left[ x^{ \star \, 0 }(\bar{\bs\sigma}); \, \bs{x}^{\star}(\bar{\bs\sigma}) \right]$ {denotes} the optimal solution to $\mc{P}(\bar{\bs\sigma})$.  Moreover, by Slater's constraint qualification in Standing Assumption \ref{ass:standing}, there exists a unique optimal dual multiplier $\lambda_{\bar{\bs\sigma}}^{\star} \in \R^n$ \cite[Section 5.2.3, p. 227]{boyd:vandenberghe}.

We now define $\left[ \bar{x}^1 \,  ; \, \ldots \,  ;  \, \bar{x}^N \right] := \bs{x}^{\star}( \bar{\bs{\sigma}} )$ and 
$\bar{\lambda} := \lambda_{\bar{\bs\sigma}}^{\star}$, 
so that at convergence (i.e., as $t \rightarrow \infty$) by \eqref{eq:dual-decomposition-bold-1}
we have
$\bar{x}^i = 
\arg\min_{y \in \mc{X}^i} J^i\left( y, \frac{1}{N} \sum_{j = 1}\bar{x}^j, \bar{\lambda} \right)$ 
for all $i \in \N[1,N]$ by Equation \eqref{eq:dual-decomposition-bold-1}, and $\frac{1}{N} \sum_{i=1}^{N} \bar{x}^i = \bs{x}^{0 \star}( \bar{\bs{\sigma}} ) \in \mathcal{S}$ by \eqref{eq:dual-decomposition-bold-2}. Thus, by Definition \ref{def:ANESC}, the pair $\left( \{\bar{x}^i\}_{i=1}^N, \, \bar{\lambda}\right)$ is an aggregative equilibrium for the game in \eqref{eq:best-response-price} with shared constraint in \eqref{eq:shared-constraint}.

\hfill $\blacksquare$

\subsection{Proof of Theorem \ref{th:eps-ANESC}}
\label{app:theorem-1}
Let, for all $i \in \N[1,N]$, $\bar{x}^i := x^{i \star}( C \bar\sigma + K \bar\lambda)$ 
as in \eqref{eq:optimal-responses}, where $\bar{\sigma} := \frac{1}{N} \sum_{i=1}^{N} \bar{x}^i$ and $\bar{\lambda}$ is fixed.
Since, for all $i$, the function $f^i$ is strongly convex, let $\bar{N} \in \N$ be such that the function $\tilde{f}^i := f^i + \tfrac{1}{2}\left\| \cdot \right\|_{ (C+C^\top)/2}^2$ is strongly convex as well, for all $i$ and for all $N \geq \bar{N}$. Thus, we define the single-valued best response
$$  \bar{x}^{i ,N}_{ \text{best} } := \underset{y \in \mc{X}^i }{\arg\min} \textstyle \, J^i\left(y, \frac{1}{N} y + \frac{1}{N} \sum_{j \neq i}^N \bar{x}^j, \bar{\lambda} \right).$$

It follows from Lemma \ref{lem:optimizer} that $\bar{x}^{i} = \left( \partial f^i\right)^{-1}\left( -C \bar{\sigma} - {K}  \bar{\lambda} \right)$, and $\bar{x}_{\textup{best}}^{i,N}=\left( \partial \tilde{f}^i \right)^{-1}\left( -C (\bar{\sigma} - \frac{1}{N} \bar{x}^{i} ) - {K}  \bar{\lambda} \right)$. 
Next, we exploit the following facts: $\left( \partial f^i\right)^{-1}$ is $(1/\ell)$-Lipschitz continuous by Lemma \ref{lem:strong-convexity-monotonicity}; $\rge ( ( \partial \tilde{f}^i )^{-1}) \subseteq \mc{X}^i \subseteq \mc{X}$ and $c_{\mathcal{X}} := \max_{x \in \mathcal{X}} \left\| x \right\| < \infty $; $\left\|   \partial f^i - \partial \tilde{f}^i \right\| \leq \tfrac{1}{N} \left\| (C + C^\top)/2 \right\| $:
\begin{align*}
& \textstyle \left\| ( \partial f^i )^{-1} - ( \partial \tilde{f}^i)^{-1}  \right\| \\
& = \textstyle \left\| \left( ( \partial f^i )^{-1} - ( \partial \tilde{f}^i)^{-1} \right) \circ ( \partial \tilde{f}^i ) \circ ( \partial \tilde{f}^i )^{-1}  \right\| \\
&\leq \textstyle \left\| ( \partial f^i )^{-1} \circ  \partial \tilde{f}^i \, - \Id \right\| \cdot \left\| ( \partial \tilde{f}^i )^{-1} \right\| \\
&\leq c_{\mathcal{X}} \textstyle \left\| ( \partial f^i )^{-1} \circ ( \partial \tilde{f}^i ) - ( \partial f^i )^{-1} \circ ( \partial f^i  ) \right\| \\
&\leq \textstyle \frac{c_{\mathcal{X}}}{\ell} \left\|   \partial \tilde{f}^i  -  \partial f^i   \right\| \\
&\leq \textstyle \frac{1}{N} \frac{c_{\mathcal{X}}}{\ell} \left\| (C+C^\top)/2 \right\|.
\end{align*}
To conclude the proof, we introduce the shorthand notation $v := -C \bar{\sigma} - {K} \bar{\lambda}$ and $v^{i} := -C (\bar{\sigma} - \frac{1}{N} \bar{x}^{i} ) - {K} \bar{\lambda}$, where $\left\| v - v^{i} \right\| = \left\| \frac{1}{N} C \bar{x}^i\right\| 
\leq \textstyle \frac{1}{N} \left\| C \right\| c_{\mathcal{X}}$.
For all $i$, we derive that
\begin{align*}
\textstyle \left\| \bar{x}^{i} - \bar{x}^{i ,N}_{ \text{best} } \right\| & = 
\textstyle \left\|  (\partial f^i)^{-1}(v) - (\partial \tilde{f}^i)^{-1}(v^i) \right\| \\
& \leq \textstyle \left\| (\partial f^i)^{-1}\left( v\right) - (\partial f^i)^{-1}\left( v^i \right) \right\| \\   
& \quad + \textstyle \left\| \left(\partial f^i\right)^{-1}\left( v^i \right) - ( \partial \tilde{f}  )^{-1}\left( v^i \right)  \right\| \\ 
& \leq c/N \,,
\end{align*}
where $c:= 3 \left\| C \right\| c_{\mathcal{X}}/{\ell}$.
\hfill $\blacksquare$

\smallskip
\begin{remark}
Theorem \ref{th:eps-ANESC} immediately implies that the difference between the optimal cost at an aggregative equilibrium and that at a Nash equilibrium vanishes in the limit of infinite population size, that is, whenever $\left( ( \bar{x}^i )_{i=1}^{N} , \, \bar{\lambda} \right)$ is an aggregative equilibrium, we have that
\begin{multline*}
\lim_{N  \rightarrow \infty}  \ \max_{i \in \N[1,N]} \textstyle \left| J^i\left( \bar{x}^i, \tfrac{1}{N} \sum_{j=1}^N \bar{x}^j, \, \bar{\lambda} \right) \right. \\
\left. -  \inf_{y \in \mathcal{X}^i} \textstyle J^i\left( y, 	\tfrac{1}{N}\left( y + \sum_{j\neq i}^{N} \bar{x}^j \right), \, \bar{\lambda} \right) \right| = 0.
\end{multline*}

In addition, if the functions $\{f^i \}_{i=1}^N$ are Lipschitz continuous, then $\{J^i \}_{i=1}^N$ in \eqref{stand:J} are Lipschitz continuous as well, hence Theorem \ref{th:eps-ANESC} implies that there exists $d \in \R_{>0}$ such that
\begin{multline*}
\max_{i \in \N[1,N]} \textstyle \left| J^i\left( \bar{x}^i, \tfrac{1}{N} \sum_{j=1}^N \bar{x}^j, \bar{\lambda} \right) \right. 
\\
\left. -  \inf_{y \in \mathcal{X}^i} \textstyle J^i\left( y, 	\tfrac{1}{N}\left( y + \sum_{j\neq i}^{N} \bar{x}^j \right), \bar{\lambda} \right) \right| \leq d/N 
\end{multline*}
for all $N \in \N$.
{\hfill $\square$}
\end{remark}

\balance

\bibliographystyle{IEEEtran}
\bibliography{library}

\end{document}